\documentclass{article}       
\usepackage{mathtools}
\usepackage{graphicx,array}
\usepackage{amssymb,amsmath,amsthm}
\usepackage{url}
\usepackage{multirow}
\usepackage{subfigure}
\usepackage{epstopdf}
\usepackage{color}
\usepackage{amsfonts}

\newcommand{\D}{\mathbb{D}}

\newcommand{\I}{\mathbb{I}}
\newcommand{\R}{\mathbb{R}}

\newcommand{\N}{\mathbb{N}}

\newtheorem{theorem}{Theorem}
\newtheorem{remark}{Remark}
\newtheorem{definition}{Definition}

\usepackage{hyperref}
\usepackage{url}

\newcommand{\splitatcommas}[1]{%
  \begingroup
  \begingroup\lccode`~=`, \lowercase{\endgroup
    \edef~{\mathchar\the\mathcode`, \penalty0 \noexpand\hspace{0pt plus 1em}}%
  }\mathcode`,="8000 #1%
  \endgroup
}

\begin{document}

\title{Trends, Directions for Further Research, and Some Open Problems of Fractional Calculus 
}


\author{Kai Diethelm%
\footnote{Faculty of Applied Natural Sciences and Humanities (FANG), University of Applied Sciences 
W\"{u}rzburg-Schweinfurt, Ignaz-Sch\"{o}n-Str.\ 11, 97421 Schweinfurt, Germany, 
kai.diethelm@fhws.de, ORCID: 0000-0002-7276-454X}\ \footnote{corresponding author}
\and Virginia Kiryakova\footnote{Institute of Mathematics and Informatics, Bulgarian Academy of Sciences, Sofia -- 1113, Bulgaria, virginia@diogenes.bg, ORCID: 0000-0002-4591-6958}
\and Yuri Luchko\footnote{Beuth Technical University of Applied Sciences Berlin, Department of  Mathematics, Physics, and Chemistry, Luxemburger Str.\ 10,  13353 Berlin, Germany, luchko@beuth-hochschule.de, ORCID: 0000-0002-9583-7377}
\and J. A. Tenreiro Machado\footnote{Institute of Engineering, Polytechnic of Porto, 
Dept.\ of Electrical Engineering, Rua Dr.\ Ant\'{o}nio Bernardino de Almeida, 431, 4249 -- 015 Porto, Portugal, 
jtm@isep.ipp.pt, ORCID: 0000-0003-4274-4879} 
\and Vasily E. Tarasov\footnote{Skobeltsyn Institute of Nuclear Physics, Lomonosov Moscow State University, Moscow 119991, Russia, tarasov@theory.sinp.msu.ru, 
ORCID: 0000-0002-4718-6274}\ \footnote{Faculty ``Information Technologies and Applied Mathematics",  
Moscow Aviation Institute (National Research University), Moscow 125993, Russia}
}


\date{09 August 2021}

\maketitle

\begin{abstract}
The area of fractional calculus (FC) has been fast developing and is presently being applied in all scientific fields. Therefore, it is of key relevance to assess the present state of development and to foresee, if possible, the future evolution, or, at least, the challenges identified in the scope of advanced research works. This paper gives a vision about the directions for further research as well as some open problems of FC. A number of topics in mathematics, numerical algorithms and physics are analyzed, giving a systematic perspective for future research.

\emph{Keywords: }Fractional calculus;  Sonine kernels; general fractional integrals and derivatives; fractional differential equations; numerical solution; fractional dynamics

\end{abstract}


\section{Introduction}\label{s1}

In 1695 Gottfried Leibniz exchanged ideas with other mathematicians about what he called Fractional Calculus (FC). FC appeared as a misname and generalized integro-differentiation seems more adequate, but the name FC remained due to historical reasons \cite{Machado2010_poster1,Machado2010_poster2,Machado5:11,Valerio:14}. The FC generalizes the standard differential calculus to non-integer orders, real or complex. This scientific tool  remained mainly in the area of mathematics  until the last two decades  (see e.g.  \cite{Kir-TMSF1996}),  when the research community recognized its superior performance for describing many natural and artificial phenomena. For a comprehensive treatment of present day knowledge in FC readers are referred to \cite{Kochubei:19a,Kochubei:19b} for mathematics, \cite{Karniadakis:19} for numerical analysis, \cite{Tarasov:19a,Tarasov:19b} physics, \cite{Petras:19} control, and \cite{Baleanu:19a,Baleanu:19b} for applications. The FC became very popular in all branches of science and an active area of development. We recall the list of 23 problems raised by Hilbert by the beginning of the 20th century and its influence upon the succeeding decades of research. While many problems were solved, efforts are still on-going for answering the remaining ones. Therefore, formulating not only possible challenges and problems, but also pointing towards future directions of research may have a significant impact.  This manuscript gives the vision of the authors in what concerns the key upcoming issues of this fast advancing area of knowledge.

\begin{figure}[h]
\centering
\includegraphics[width=1.0\linewidth]{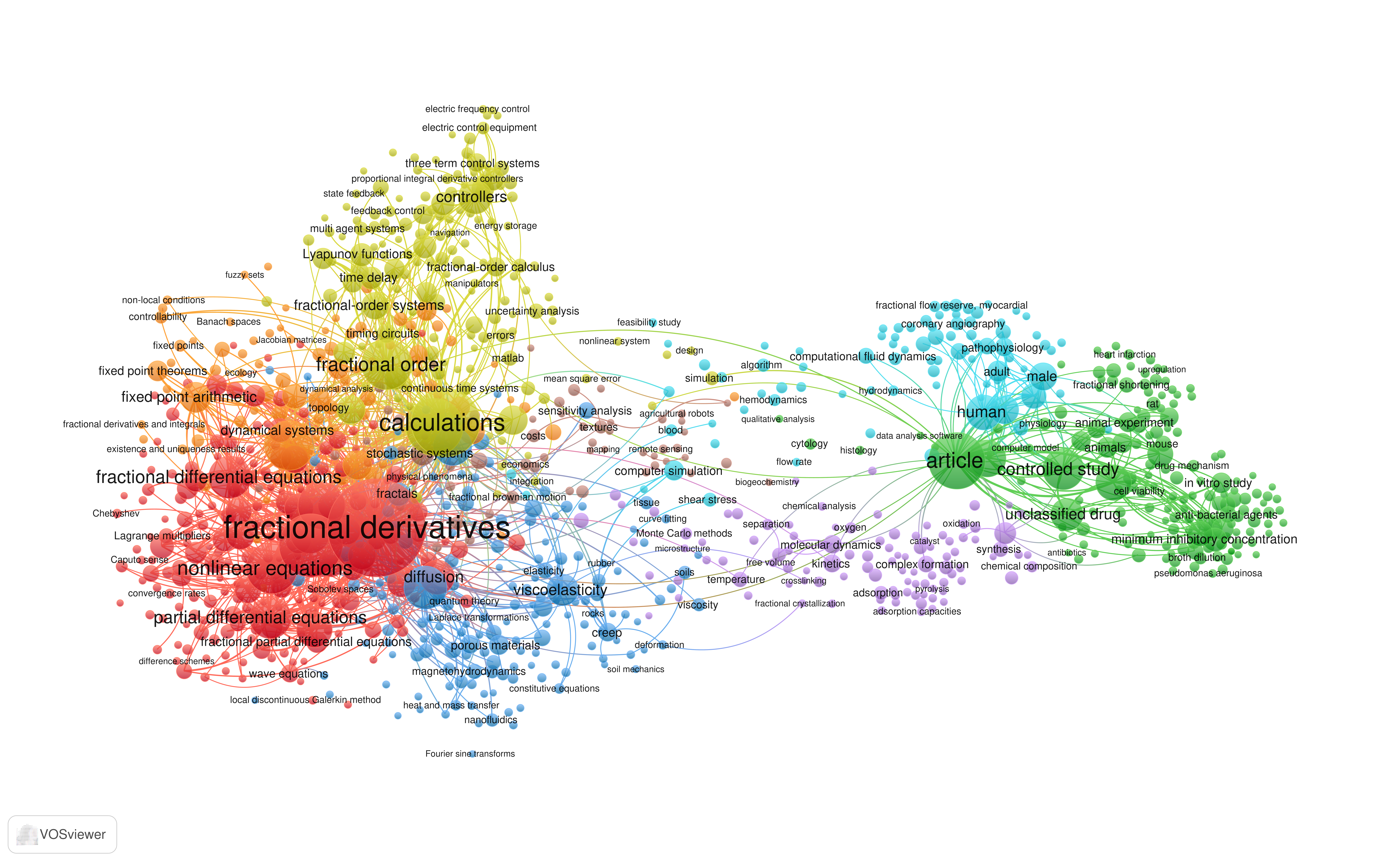}
\caption{Network plot using  \textsf{Co-occurrence}, with options \textsf{Index keywords}, \textsf{Fractional counting}.}
\label{fig:Co-occurrence_Index_keywords_Fractional_counting}
\end{figure}
\begin{figure}[h]
\centering
\includegraphics[width=1.0\linewidth]{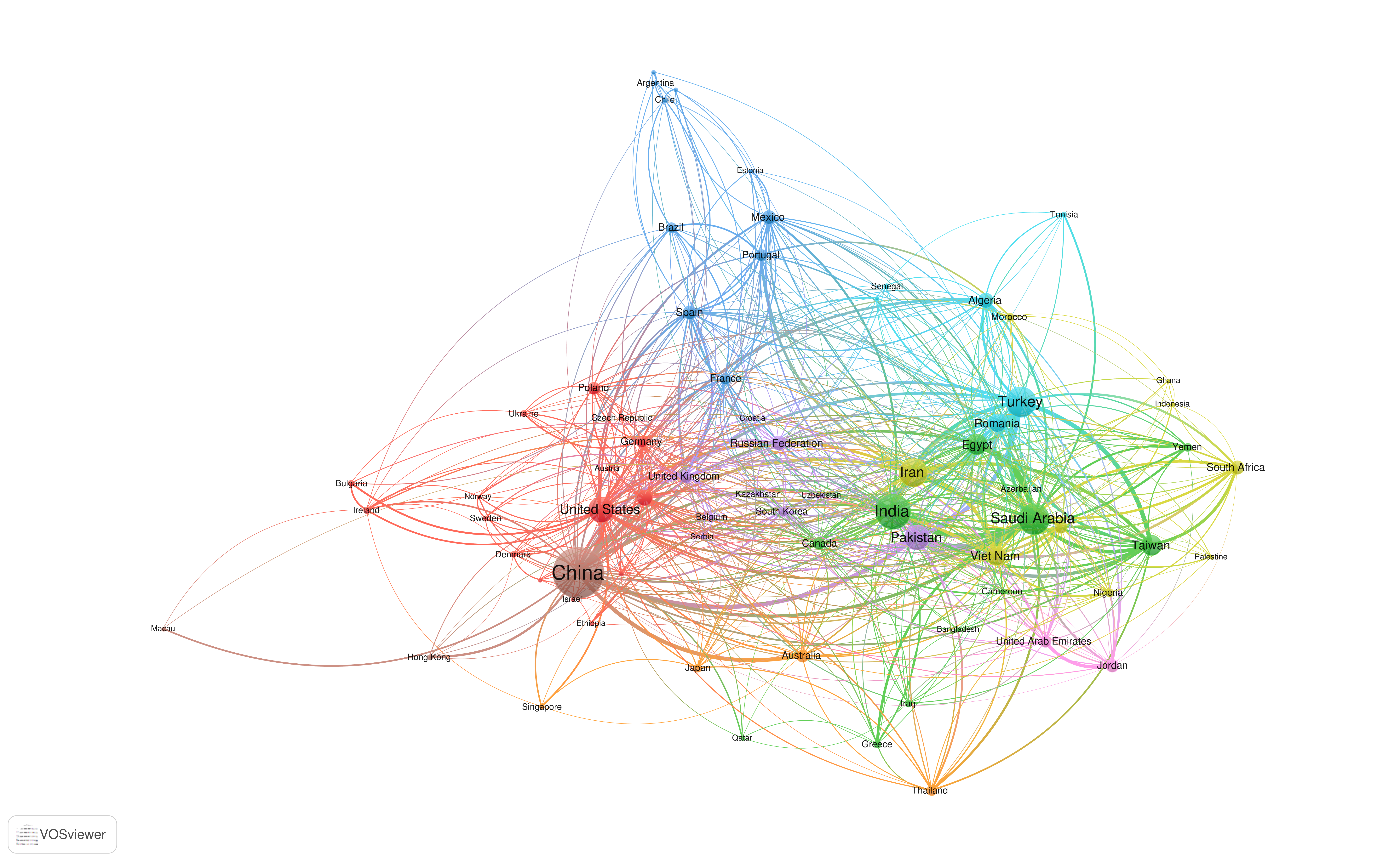}
\caption{Network plot using \textsf{Co-authorship}, with options \textsf{Countries}, \textsf{Fractional counting}.}
\label{fig:Co-authorship_Countries_Fractional_counting}
\end{figure}
\begin{figure}[h]
\centering
\includegraphics[width=1.0\linewidth]{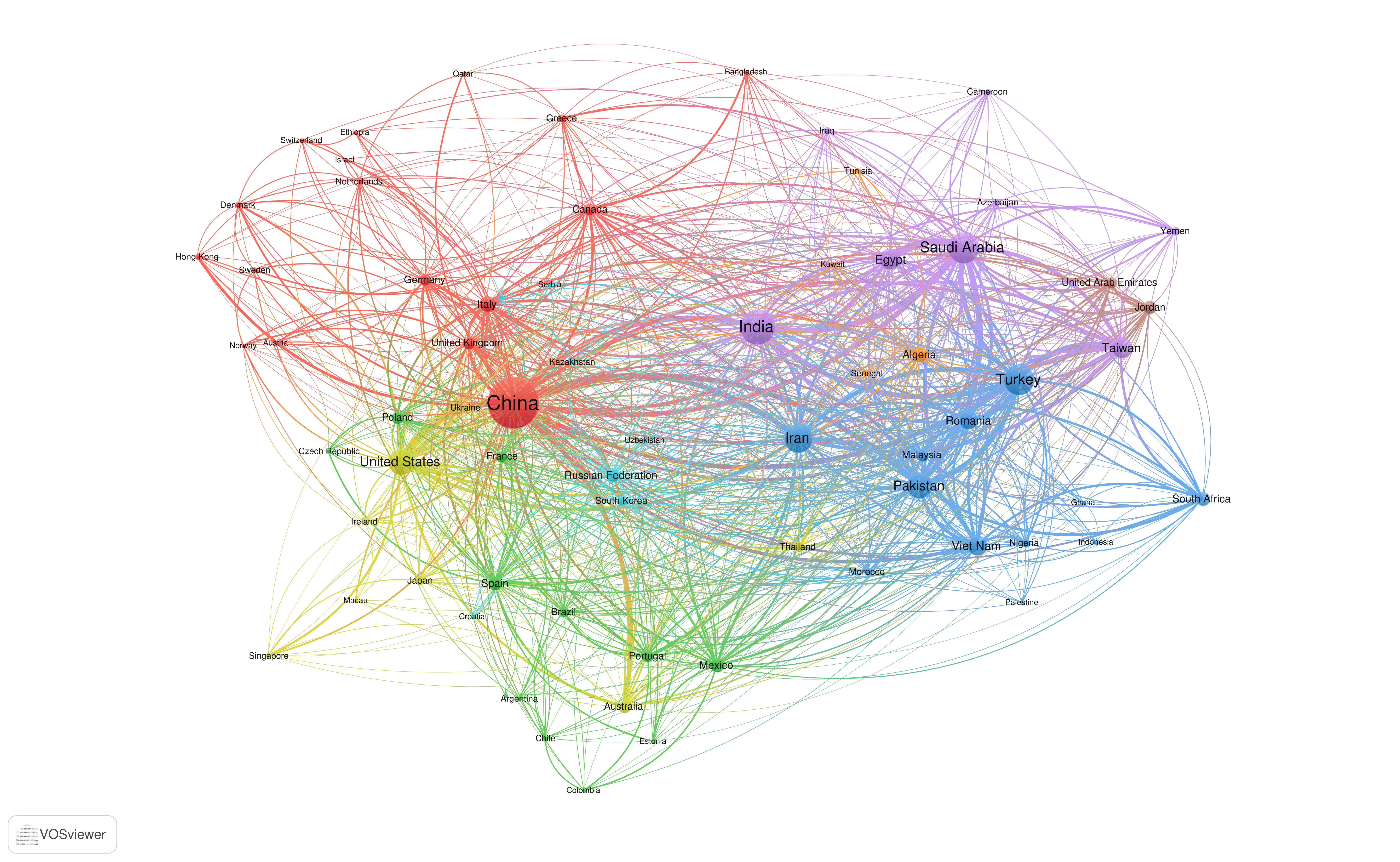}
\caption{Network plot using \textsf{Bibliographic coupling}, with options \textsf{Countries}, \textsf{Fractional counting}.}
\label{fig:Bibliographic_coupling_Countries}
\end{figure}

We can estimate of the present day state of FC using present day publicly available information, just to remind that until 1974 there were only 1 book  and 1 conference proceedings, devoted to FC as a topic, while by 2018 the FC books were estimated to more than 240 (see \cite{Machado-Kir-FCAA2017} and Table 1 in \cite{Machado-Kir-HFCA-Ch1}), and the published FC articles hardly to be  count as thousands. For that purpose we selected the program VOSviewer \cite{vanEck:2009,Waltman:2010,vanEck:2014,Rodriguez:2016,vanEck:2017,MoralMuoz:2020} as the software tool for processing and clustering bibliographic information. We tackled data available at Scopus database and collected papers published during year 2020. Moreover, 8 search keywords where adopted, namely the set  \{Fractional calculus, Fractional derivative, Fractional integration, Fractional dynamics, Mittag-Leffler, Derivative of non-integer order, Integral of non-integer order, Derivative of complex order, Integral of complex order\} that led to 6,589 records. The VOSviewer allows a multitude of possible perspectives of bibliographic analysis, but, for the sake of parsimony, we focus on 3 types of scrutiny based on network plots:


\begin{itemize}
\item \textsf{Co-occurrence}, with options \textsf{Index keywords}, \textsf{Fractional counting}, \textsf{Minimum number of occurrence of a keyword}: 5, that reveals 1,116 keywords, as portrayed in Fig. \ref{fig:Co-occurrence_Index_keywords_Fractional_counting} 
\item \textsf{Co-authorship}, with options \textsf{Countries}, \textsf{Fractional counting}, \textsf{Minimum number of documents of a country}: 5, and \textsf{Minimum number of citations of a country}: 2, that includes 71 countries, represented in Fig. \ref{fig:Co-authorship_Countries_Fractional_counting}
\item \textsf{Bibliographic coupling}, with options Countries, \textsf{Fractional counting}, \textsf{Minimum number of documents of a country}: 5, \textsf{Minimum number of citations of a country}: 2, that leads to 71 countries, shown in Fig. \ref{fig:Bibliographic_coupling_Countries} 
\end{itemize}

Figure \ref{fig:Co-occurrence_Index_keywords_Fractional_counting} shows that  FC is presently applied in all fields of science. Additionally, we observe two main clusters: one that includes the areas of mathematics, physics, engineering and economy, and a more recent one that covers medicine, biology and genetics. Figures \ref{fig:Co-authorship_Countries_Fractional_counting} and \ref{fig:Bibliographic_coupling_Countries}  portrait also two main clusters. One major cluster corresponding to China, USA, most European countries, Russia, south America and several countries of Asia, such as, Japan, South Korea, and Australia. A second cluster is formed by more recent contributions mostly from  Middle East and south Asia, such as India, Iran, Turkey, Saudi Arabia and Pakistan.

The emerging clusters are clearly correlated with the world state of FC. On one hand, it is well known the applicability of FC is describing complex phenomena including non-locality and memory effects, which accounts for the results in Fig. \ref{fig:Co-occurrence_Index_keywords_Fractional_counting} . On the other hand, the portraits in Figs. \ref{fig:Co-authorship_Countries_Fractional_counting} and \ref{fig:Bibliographic_coupling_Countries} are also natural, since we have an historical development of FC  in the countries located in the first  cluster, while we observe countries joining more recently  the area of FC that form the right cluster.

Given this state of affairs, we outline our vision for a future solid progress of FC in the next sections. Section \ref{s2}  discusses important analytical aspects that need further work. Section \ref{s3} addresses  numerical problems in the scope of present day computational era.  Section \ref{s4} analyses problems in physics and the application of the FC tools. Finally, Section \ref{s5} outlines the main conclusions. Due to space limitations, these sections cover a limited set of topics. Nonetheless, the present work gives a synthetic view and stimulates not only a sound scientific debate, but also gives room for a future continuation.


\section{Analytical Aspects} \label{s2}


The fractional integrals and derivatives of arbitrary order $\alpha \in \mathbb{R}$  or even $\alpha \in \mathbb{C}$ can be roughly characterized as an interpolation of the $n$-th order derivatives and the $n$-fold definite integrals ($n\in \mathbb{N}$). Evidently, without posing some  additional conditions, this interpolation is not unique and thus, several different definitions of the FC operators have been suggested starting from the very origins of this theory. In particular, the Riemann-Liouville integral and derivative, the Gr\"{u}nwald-Letnikov derivative, the Marchaud derivative, the Weyl integral and derivative of periodic functions, the Erd\'{e}lyi-Kober integral and derivative, the Hadamard integral and derivative, the Riesz and the Feller potentials and fractional derivatives, and the Dzherbashian-Caputo fractional derivative have been introduced and studied to mention only some of the most used FC fractional integrals and derivatives. However, these operators were defined neither artificially nor arbitrary. All of them possess a strong mathematical motivation and their definitions and properties were exactly adjusted both to the spaces of functions and their domains and to the problems that these operators were useful for. For a detailed discussion of the classical FC operators including their history and applications we refer the readers to the ``bible of FC'' - the encyclopedic monograph \cite{SKM} and to the recent books \cite{Kochubei:19a,BDST2016,Di2010}. Later on, some other fractional integrals and derivatives were introduced in form of integral and integro-differential operators with the suitable special functions in the kernels, see e.g. \cite{VKir-Book,Kir-HFCA-Ch4}.  These operators proved to have also various useful applications in several mathematical topics.    

Within the last few years, the situation in FC essentially changed. Several new definitions of the fractional  derivatives were suggested and, unfortunately, some of them cannot be called neither fractional nor derivatives, see, e.g.,  \cite{DGGS}, \cite{Han}. Thus, the  questions like ``What are the fractional integrals and derivatives?'', ``What are their decisive mathematical properties?'', ``What fractional operators make sense in applications and why?'' became extremely important and acute.

The aim of this section is in providing some partial answers to the questions mentioned above (see \cite{Luc21} for further contributions to this topic). We start with a discussion of the characteristic properties that the one-parameter FC operators should possess. Then we consider an example of  realization of this schema in the case of the fractional integrals and derivatives acting on the functions defined on a finite interval. Finally, we discuss a class of the general fractional derivatives with the Sonine kernels defined for the functions on the real positive semi-axis and formulate some open problems for further research.

\subsection{Desiderata for the FC operators}
\label{sec2.2}

In this subsection, we closely follow the ideas recently presented  in \cite{HL}. For the sake of brevity, we restrict ourselves to the case of the one-parameter fractional integrals and derivatives of a real order $\alpha \ge 0$. The properties formulated below should be valid on some suitable spaces of functions that we do not specify here (see \cite{HL} for details). 

The  families of the operators  $I^\alpha,\ \alpha \ge 0$ and $D^\alpha,\ \alpha \ge 0$  can be interpreted as  fractional integrals and fractional derivatives, respectively,  if they satisfy the following six properties (called desiderata in \cite{HL}):

\vspace{0.2cm}

\noindent
HL1) The operators $I^\alpha$  and $D^\alpha$ are linear operators on the suitable linear spaces of functions.

\vspace{0.2cm}

\noindent
HL2) For the operators  $I^\alpha$, 
the index law (semigroup property)
\begin{equation}
\label{s2_1}
I^{\alpha_1} \circ I^{\alpha_2}=I^{\alpha_1+\alpha_2}, \ \ \alpha_1, \alpha_2 \ge 0
\end{equation}
holds true, where
$\circ$ denotes composition of operators.

\vspace{0.2cm}

\noindent
HL3) The operator $D^\alpha$ is left inverse  to the operator $I^\alpha$:
\begin{equation}
\label{s2_2}
D^{\alpha} \circ I^{\alpha}=\mbox{Id}, \ \ \alpha \ge 0.
\end{equation}

\vspace{0.2cm}

\noindent
HL4) The following limits exist in some sense
\begin{equation}
\label{s2_3}
\lim_{\alpha\to 0}D^{\alpha} = D^0,\ \ \lim_{\alpha\to 1}D^{\alpha} = D^1
\end{equation}
and define linear maps. 

\vspace{0.2cm}

\noindent
HL5) The map $D^0$ is the identity operator.

\vspace{0.2cm}

\noindent
HL6) The map $D^1$ is the first order derivative. 

\vspace{0.2cm}

We mention once again that the precise formulations of the properties (desiderata) HL1)-HL6) presented above can be found in \cite{HL}. In particular, in \cite{HL}, a description of the spaces of functions, where these properties should hold valid, is provided. Some examples of realization of the schema  specified by the properties HL1)-HL6) were discussed in the recent publications \cite{HK,KH,KH1}.

\subsection{One-parameter FC operators on a final interval}
\label{sec2.3}

This subsection is devoted to a description of a one-parameter family of the fractional integrals defined for the functions on a finite interval and several families of the corresponding fractional derivatives in the sense of the desiderata presented in the previous subsection. For the details and the proofs of the theorems we refer to the recent publication \cite{Luc20}.

We start with an important result regarding a characterization of a one-parameter family of the fractional integrals defined for the functions on a finite interval.

\begin{theorem}[\cite{Cart78}]
\label{t_1}
Let $E$ be the space $L_p(0,1),\, 1\le p <+\infty$, or $C[0,1]$. Then there is precisely one family $I^\alpha,\, \alpha >0$ of operators on $E$ satisfying the following conditions:

\vspace{0.2cm}

\noindent
CM1)  $(I^1\, f)(x) = \int_0^x f(t)\, dt, \ f\in E$ (interpolation condition),
\vspace{0.1cm}

\noindent
CM2) $(I^\alpha \, I^\beta\, f)(x) = (I^{\alpha + \beta}\, f)(x), \ \alpha,\beta > 0, \  f\in E$ (index law),
\vspace{0.1cm}

\noindent
CM3) $\alpha \to I^\alpha$ is a continuous map of $(0,\, +\infty)$ into ${\mathcal L}(E)$ for some Hausdorff topology on ${\mathcal  L}(E)$, weaker than the norm topology (continuity),

\noindent
CM4) $f\in E$ and $f(x)\ge 0$ (a.e. for $E=L_p(0,1)$) $\Rightarrow$ $(I^\alpha\, f)(x) \ge 0$ (a.e. for $E=L_p(0,1)$) for all $\alpha > 0$ (non-negativity).

\vspace{0.2cm}

\noindent
That family is given by the Riemann-Liouville fractional integrals with $\alpha >0$:
\begin{equation}
\label{RLI}
(I^\alpha\, f)(x) = \begin{cases}
\frac{1}{\Gamma(\alpha)} \int_0^x (x-t)^{\alpha -1}\, f(t)\, dt, & \alpha >0,\\
f(x), & \alpha = 0.
\end{cases}
\end{equation}
\end{theorem}

As we see, Theorem \ref{t_1} does not cover the case $\alpha =0$.  However, it is well known (Theorem 2.6 in \cite{SKM})
that the Riemann-Liouville fractional integrals $I^\alpha,\, \alpha \ge 0$ given by \eqref{RLI} form a semigroup in $L_p(0,\ 1),\ p\ge 1$, which is strongly continuous for all $\alpha \ge 0$, i.e., the relation
\begin{equation}
\label{cont}
\lim_{\alpha \to \alpha_0} \| I^\alpha f - I^{\alpha_0} f\|_{L_p(0,1)} = 0
\end{equation}
holds valid for any $\alpha_0,\ \alpha_0 \ge 0$ and for any $f\in L_p(0,1)$.  Thus, we can uniquely extend the family of the Riemann-Liouville fractional integrals defined for $\alpha >0$ to a family defined for $\alpha \ge 0$ by setting $I^0 = \mbox{Id}$. Evidently, this extended family of operators satisfies the properties CM1)-CM4) for all $\alpha \ge 0$.

The conditions formulated in Theorem \ref{t_1} are much stronger than those  presented in the previous subsection. However, they are very natural and useful for applications and thus, in a certain sense, the Riemann-Liouville fractional integrals are the only ``right" one-parameter family of the fractional integrals defined for the functions on a finite interval.

In what follows, without loss of generality, we restrict ourselves to the operators defined on the interval $[0,\, 1]$.

In order to introduce the suitable fractional derivatives, the following result regarding solvability of the Abel integral equation is employed:

\begin{theorem}[\cite{SKM}]
\label{t_A}
The Abel integral equation 
\begin{equation}
\label{Abel}
(I^\alpha\, \phi)(x) = \frac{1}{\Gamma(\alpha)} \int_0^x (x-t)^{\alpha -1}\, \phi(t)\, dt\ = f(x),\ \,
0<\alpha<1,\ x\in[0,\, 1]
\end{equation}
possesses a unique solution in $L_1(0,1)$ if and only if
\begin{equation}
\label{Acond}
I^{1-\alpha}\, f\in \mbox{AC}([0,1]) \ \ \mbox{and} \ \ (I^{1-\alpha}\, f)(0) = 0,
\end{equation}
where the notation $\mbox{AC}([0,1])$ stands for the space of functions that are absolutely continuous on the interval $[0,\, 1]$:
\begin{equation}
\label{absc}
f\in \mbox{AC}([0,1]) \, \Leftrightarrow \, \exists \phi \in L_1(0,1):\, f(x) = f(0) + \int_0^x \phi(t)\, dt,\ x\in [0,\, 1].
\end{equation}
If these conditions are satisfied, the solution is given by the formula 
\begin{equation}
\label{Asol}
\phi(x)  = \frac{d}{dx}\, (I^{1-\alpha}\, f)(x) = \frac{d}{dx}\, \frac{1}{\Gamma(1-\alpha)} \int_0^x (x-t)^{-\alpha}\, f(t)\, dt,\ \, 0<\alpha<1.
\end{equation}
\end{theorem}

Using  the representation \eqref{absc}, we define a (weak) derivative of a function $f\in \mbox{AC}([0,1])$:
\begin{equation}
\label{absc_d}
f(x) = f(0) + \int_0^x \phi(t)\, dt,\ x\in [0,\, 1] \ \Rightarrow \ \frac{df}{dx} := \phi \in L_1(0,1).
\end{equation}
In the further discussions, we interpret the first order derivative of the absolutely continuous functions in the sense of the formula  \eqref{absc_d}.

According to the property HL3) from the desiderata, a fractional derivative for the functions defined on a finite interval is introduced as a linear operator $D^\alpha,\ \alpha \ge 0$ left inverse  to the Riemann-Liouville fractional integral $I^\alpha$, i.e., as an operator that satisfies the relation (1st Fundamental Theorem of FC)


\begin{equation}
\label{ftFC}
(D^\alpha\, I^\alpha\, \phi)(x) = \phi(x),\ x\in [0,\, 1]
\end{equation}
on appropriate nontrivial space of functions.

It turns out  that there exist infinitely many different families of the fractional derivatives in the sense of relation \eqref{ftFC}. In the rest of this subsection, we discuss some known and new families of the fractional derivatives on the interval $[0,\, 1]$ (see \cite{Luc20} for details). 

\begin{remark}
\label{r_1}
The formula \eqref{ftFC} and the relation $(I^0\, f)(x) = f(x)$ uniquely define the fractional derivative of order $\alpha =0$ as the identity operator: $(D^0\, f)(x) = f(x)$. In what follows, we mainly restrict ourselves  to the case $0< \alpha \le 1$ and to the space  $L_1(0,1)$ and its subspaces (a similar theory can be developed for, say, $L_p(0,1),\ 1<p<+\infty$ and its subspaces).
\end{remark}

The formula \eqref{ftFC} can be  rewritten in equivalent form of two equations:
\begin{equation}
\label{ftFC_1}
(I^\alpha\ \phi)(x) = f(x),\ \ (D^\alpha\, f)(x) = \phi(x),\ \,  x\in [0,\, 1].
\end{equation}
As we see, the second of equations \eqref{ftFC_1} defines the fractional derivative $D^\alpha$ of a function $f$ as the solution $\phi$ of the Abel integral equation with the right-hand side $f$. Now we recall Theorem \ref{t_A} and reformulate it as follows:

\begin{theorem}
\label{t_A-1}
On the space of functions $I^\alpha(L_1(0,1))$, the unique fractional derivative $D^\alpha$ of order $\alpha$, $0< \alpha < 1$ is the Riemann-Liouville fractional derivative given by the formula
\vskip -12pt
\begin{equation}
\label{RLD}
(D^\alpha_{RL}\, f)(x)  = \frac{d}{dx}\, (I^{1-\alpha}\, f)(x) = \frac{d}{dx}\, \frac{1}{\Gamma(1-\alpha)} \int_0^x (x-t)^{-\alpha}\, f(t)\, dt.
\end{equation}
\end{theorem}

Evidently, the formula \eqref{RLD} makes sense for a space of functions larger than $I^\alpha(L_1(0,1))$, namely, for the space $\{ f:\, I^{1-\alpha}\, f \in \mbox{AC}([0,\, 1])\}$.

Moreover, the 1st Fundamental Theorem of FC (relation \eqref{ftFC}) for the Riemann-Liouville fractional derivative is valid on the whole space $L_1(0,1)$, i.e., the formula
\begin{equation}
\label{ftFC_RL}
(D^\alpha_{RL}\, I^\alpha\, f)(x) = f(x),\ x\in [0,\, 1],\ f \in L_1(0,1)
\end{equation}
holds true (\cite{SKM}).

For the Riemann-Liouville fractional derivatives, the properties HL4)-HL6) from the desiderata formulated in the previous subsection are valid on the suitable spaces of functions (see, e.g., \cite{SKM}). Thus, the families of the Riemann-Liouville fractional integrals and derivatives are the  fractional integrals and derivatives in the sense of the desiderata. 

However, there exist infinitely many other one-parameter families of the fractional derivatives associated to the Riemann-Liouville fractional integrals. They are introduced in the rest of this subsection (see \cite{Luc20} for details). 

On the space $\mbox{AC}([0,\, 1])$, the Caputo fractional derivative $D^\alpha_C$ of order $\alpha$, $0 < \alpha \le 1$ is defined as follows:
\begin{equation}
\label{CD}
(D^\alpha_{C}\, f)(x) = (I^{1-\alpha} \frac{d}{dx}f)(x).
\end{equation}

As in the case of the Riemann-Liouville fractional derivative, the 1st Fundamental Theorem of FC (relation \eqref{ftFC}) for the Caputo fractional derivative is valid on even larger space of functions
\begin{equation}
\label{X2}
X_{FT} = \left\{ f:\, I^{\alpha}f\in \mbox{AC}([0,\, 1])\ \mbox{and}\ (I^{\alpha}f)(0) = 0\right\},
\end{equation}
i.e., the formula
\begin{equation}
\label{ftFC_C}
(D^\alpha_{C}\, I^\alpha\, f)(x) = f(x),\ \, x\in [0,\, 1],\ f \in X_{FT}
\end{equation}
holds true (\cite{Luc20}). It is worth mentioning that the space $X_{FT}$ can be also characterized as follows (Theorem 2.3 in \cite{SKM}):
\begin{equation}
\label{X2_1}
X_{FT}= I^{1-\alpha}(L_1(0,1))\, (\forall f \in X_{FT}\, \exists \phi \in L_1(0,1):\, f(x) = (I^{1-\alpha}\, \phi)(x)).
\end{equation}

Another family of the fractional derivatives associated to the Riemann-Liouville fractional integral is the  Hilfer fractional derivatives of order  $\alpha,\ 0<\alpha \le 1$ and type $\gamma_1$, $0\le \gamma_1 \le 1-\alpha$:
\begin{equation}
\label{HD}
(D^{\alpha,\gamma_1}_H\, f)(x)  = (I^{\gamma_1}\, \frac{d}{dx}\, I^{1-\alpha-\gamma_1}\, f)(x).
\end{equation}
This derivative is well defined on the space 
$\left\{ f:\, I^{1 - \alpha -\gamma_1}\, f \in \mbox{AC}([0,\, 1])\right\}$. 

The 1st Fundamental Theorem of FC (relation \eqref{ftFC}) for the Hilfer fractional derivative is valid on the  space $X_{FT}$ defined by \eqref{X2} (\cite{Luc20}):
\begin{equation}
\label{ftFC_H}
(D^{\alpha,\gamma_1}_H\, I^\alpha\, f)(x) = f(x),\ x\in [0,\, 1],\ f \in X_{FT},\ 0 < \alpha \le 1,\ 0\le \gamma_1 \le 1-\alpha.
\end{equation}

For every type $\gamma_1$, $0 \le \gamma_1 \le 1- \alpha$, the Hilfer fractional derivatives $D^{\alpha,\gamma_1}_H$ of order $\alpha,\ 0< \alpha \le 1$ form a one-parameter family of the fractional derivatives in the sense of the 1st Fundamental Theorem of FC (relation \eqref{ftFC}).  For $\gamma_1 = 0$, this family coincides with the Riemann-Liouville fractional derivatives, while for $\gamma_1 = 1-\alpha$ we get the Caputo fractional derivatives.

The construction of the Hilfer fractional derivative can be extended to the case of $n$ compositions of the first order derivatives and appropriate Riemann-Liouville fractional integrals.  As a result, we get infinitely many different families of the one-parameter fractional derivatives that were called in \cite{Luc20} the $n$th level fractional derivatives of order $\alpha,\ 0<\alpha \le 1$ and type $\gamma = (\gamma_1,\ \gamma_2,\dots,\gamma_n)$.

To define these fractional derivatives, we introduce the parameters $ \splitatcommas{  \gamma_1, \gamma_2, \dots, \gamma_n \in \R} $ that satisfy the conditions
\begin{equation}
\label{gamma}
0\le \gamma_k \ \mbox{and}\  \alpha + s_k \le k, \ \  k=1,2,\dots,n
\end{equation}
with $s_k$ in form
\begin{equation}
\label{not}
s_k:= \sum_{i=1}^k\, \gamma_i,\ k= 1,2,\dots,n.
\end{equation}

The $n$th level  fractional derivative of order  $\alpha,\ 0<\alpha \le 1$ and type $\gamma = (\gamma_1,\ \gamma_2,\dots,\gamma_n)$ 
is defined as follows (\cite{Luc20}:
\begin{equation}
\label{nLD}
(D^{\alpha,(\gamma)}_{nL}\, f)(x)  = \left(\prod_{k=1}^n (I^{\gamma_k}\, \frac{d}{dx})\right)\, (I^{n-\alpha-s_n}\, f)(x).
\end{equation}
The right-hand side of \eqref{nLD} is well defined for the functions from the space 
$\{ f: \left(\prod_{k=i}^n (I^{\gamma_k}\, \frac{d}{dx})\right)\, I^{n-\alpha-s_n}\, f \in \mbox{AC}([0,\, 1]),\ \ i=2,\dots n+1 \}$. 

The $n$th level  fractional derivative of order  $\alpha,\ 0<\alpha \le 1$ and any type $\gamma = (\gamma_1,\ \gamma_2,\dots,\gamma_n)$ is associated to the Riemann-Liouville fractional integral via the 1st Fundamental Theorem of FC.

\begin{theorem}[\cite{Luc20}]
\label{t_ft_Ln}
The $n$th level  fractional derivative is a left-inverse operator to the Riemann-Lioville fractional integral on the space $X_{FT}$ defined by \eqref{X2}:
\begin{equation}
\label{ftFC_nL}
(D^{\alpha,(\gamma)}_{nL}\, I^\alpha\, f)(x) = f(x),\ f\in X_{FT},\  x\in [0,\, 1].
\end{equation}
\end{theorem}

\begin{remark}
In  \cite{Dzherbashian2020}, the Cauchy problems for the fractional differential equations with the operators similar to the $n$th level fractional derivatives $D^{\alpha,(\gamma)}_{nL}$ (in other notations and with other restrictions on the parameters) were treated. However, the connection to the Riemann-Liouville fractional integrals in form of Theorem \ref{t_ft_Ln} was not discussed in \cite{Dzherbashian2020}.
\end{remark}

For further properties of the $n$th level  fractional derivatives and the fractional relaxation equations with these derivatives we refer the interested reader to \cite{Luc20} and \cite{Luc20a}.

\subsection{General FC operators with the Sonine kernels}
\label{sec2.4}

The previous subsections were devoted to the one-parameter families of the fractional integrals and derivatives. Without any doubts, these FC operators can be interpreted as direct generalizations of the $n$-folds integrals and integer order derivatives and thus they are objects important and useful both for analytical and numerical treatment and for applications. However, within the last years, some more general fractional operators were introduced and employed as models in several applications including the linear viscoelastisity and anomalous diffusion and wave propagation. In particular, much attention was given to the multi-term fractional differential operators that are the finite sums of the one-parameter fractional derivatives with different orders and their generalizations in form of the distributed order fractional derivatives that can be interpreted as sums of infinitely many fractional derivatives with different orders that belong to a certain interval (usually, to the intervals $[0,1]$, or $[1,2]$, or $[0,2]$). These more general fractional derivatives are determined either by finite (multi-term case) or infinite (distributed order case) many orders of the fractional derivatives and thus they require a different procedure for their characterization compared to the schema we presented in the previous subsections. 

One of the possible approaches for a description of both the single-term, the multi-term, and the distributed order fractional derivatives is in the framework of the so called general fractional calculus suggested in \cite{Koch11}. 

In this subsection, we present some results derived in \cite{Koch11} and in the recent publications \cite{Luc21a}-\cite{Luc21c} regarding  the general fractional derivatives (GFD) and the general fractional integrals (GFI). We start with the case of the operators with the ``generalized order'' restricted to the interval $(0,\ 1)$. The GFDs  of the Riemann-Liouville and Caputo types, respectively, are defined as the following integro-differential operators of the convolution type:
\begin{equation}
\label{FDR-L}
(\D_{(k)} f) (t) = \frac{d}{dt}\int_0^t k(t-\tau)f(\tau)\, d\tau,
\end{equation}
\begin{equation}
\label{FDC}
(_*\D_{(k)} f) (t) = \int_0^t k(t-\tau)f^\prime(\tau)\, d\tau,
\end{equation}
where $k$ is a nonnegative locally integrable function. For an absolutely
continuous function $f$ satisfying $f^\prime \in L^1_{loc}(\R_+)$,
the relation 
\begin{equation}
\label{rel}
(_*\D_{(k)} f) (t)\ = \ (\D_{(k)} (f-f(0))) (t) \ = \ (\D_{(k)} f) (t) - k(t)f(0)
\end{equation}
between the Caputo and Riemann-Liouville types of GFDs holds true.
In \cite{Koch11}, the  Caputo type GFD was
introduced in form of the right-hand side of \eqref{rel} that is well defined for a lager class of functions (in particular, for absolutely continuous functions) compared to 
the definition \eqref{FDC}
that requires the inclusion $f^\prime
\in L^1_{loc}(\R_+)$. In what follows, we mainly address the GFD ot the Caputo type in the sence of the right-hand side of the formula \eqref{rel}.

The Riemann-Liouville and Caputo  fractional derivatives defined by \eqref{RLD} and \eqref{CD}, respectively,  are particular cases of the GFDs \eqref{FDR-L} and \eqref{FDC} with the kernel
\begin{equation}
\label{single}
k(t) = h_{1-\alpha}(t),\ 0 <\alpha <1
\end{equation}
with the power function $h_\beta(t):= t^{\beta -1}/\Gamma(\beta),\ t>0,\ \beta >0$.
 Other important
particular cases of \eqref{FDR-L} and \eqref{FDC} are the multi-term fractional derivatives and the fractional derivatives of the
distributed order. They are generated by \eqref{FDR-L} and \eqref{FDC} with the kernels 
\begin{equation}
\label{multi}
k(t) = \sum_{k=1}^n a_k\, h_{1-\alpha_k}(t),
\ \  0 < \alpha_1 <\dots < \alpha_n < 1,\ a_k\in \R,\ k=1,\dots,n
\end{equation}
and
\begin{equation}
\label{distr}
k(t) = \int_0^1 h_{1-\alpha}(t)\, d\rho(\alpha),
\end{equation}
respectively, where $\rho$ is a Borel measure on $[0,\, 1]$.

In \cite{Koch11}, the case of the GFDs \eqref{FDR-L} and \eqref{rel} with the kernels that satisfy the following properties was considered:

\vspace{0.2cm}

\noindent
K1) The Laplace transform $\tilde k$ of $k$,
$$
\tilde k(p) = ({\mathcal L}\, k)(p)\ =\ \int_0^\infty k(t)\, e^{-pt}\, dt
$$
exists for all $p>0$,

\vspace{0.1cm}

\noindent
K2) $\tilde k(p)$ is a Stieltjes function,

\vspace{0.1cm}

\noindent
K3) $\tilde k(p) \to 0$ and $p \tilde k(p) \to \infty$ as $p \to \infty$,

\vspace{0.1cm}

\noindent
K4)  $\tilde k(p) \to \infty$ and $p \tilde k(p) \to 0$ as $p \to 0$.

\vspace{0.2cm}

In what follows, we denote the set of the kernels that satisfy the conditions K1)-K4) by $\mathcal{K}$.
As shown in \cite{Koch11}, for each $k \in \mathcal{K}$, there 
exists a completely monotone function $\kappa$ such that the Sonine condition holds true:
\begin{equation}
\label{kappa}
(k*\kappa)(t) = \int_0^t k(t-\tau) \kappa(\tau)\, d\tau \ = \ \{1\}, \ t>0,
\end{equation}
where by $\{1\}$ we denote the function that takes the constant value $1$ for $t>0$. 
The kernels that satisfy the Sonine condition \eqref{kappa} are called the Sonine kernels. The set of all Sonine kernels is denoted by $\mathcal{S}$ ($\mathcal{K} \subset \mathcal{S}$). For the properties of the Sonine kernels and their examples we refer to \cite{Han,Luc21a,Luc21b,Luc21c,Sam,Son} and other related publications.

A GFI with the kernel $\kappa$ associated to the kernel $k$ of the GFD by means of the relation \eqref{kappa} is introduced as follows:
\begin{equation}
\label{FI}
(\I_{(k)} f) (t) = \int_0^t \kappa(t-\tau)f(\tau)\, d\tau.
\end{equation}

The notions of GFD and GFI are justified by the following relations (\cite{Koch11}): 
\begin{equation}
\label{FDI}
( _*\D_{(k)} \, \I_{(k)}\, f) (t) = f(t),
\end{equation}
for any locally bounded measurable function $f$ on $\R_+$ and 
\begin{equation}
\label{FDI_2}
( \I_{(k)}\,  _*\D_{(k)}\, f) (t) = f(t) - f(0)
\end{equation}
for any absolutely continuous function $f$ on $\R_+$.

The formula \eqref{FDI} and the relation \eqref{rel} between the GFDs of the Caputo and the Riemann-Liouville types lead to the identity
\begin{equation}
\label{FDIRL}
(\D_{(k)}\, \I_{(k)} \, f) (t) = f(t),
\end{equation}
i.e., the Riemann-Liouville GFD is also a left inverse operator to the GFI defined by \eqref{FI}.

As mentioned in \cite{Han,Sam}, the functions that satisfy the Sonine condition \eqref{kappa} cannot be continuous at the point $t=0$ and thus the ``new fractional derivatives'' with the continuous kernels introduced recently in the FC literature do not belong to the class of the GFDs that we discuss in this subsection. 


Recently, in the papers \cite{Luc21a}-\cite{Luc21c}, another important class of the Sonine kernels was treated. \begin{definition}[\cite{Luc21a}]
\label{dd2}
Let $\kappa,\, k \in C_{-1,0}(0,+\infty)$, where
\begin{equation}
\label{subspace}
 C_{\alpha,\beta}(0,+\infty) \, = \, \{f:\ f(t) = t^{p}f_1(t),\ t>0,\ \alpha < p < \beta,\ f_1\in C[0,+\infty)\},
\end{equation} 
be a pair of the Sonine kernels, i.e., the Sonine condition \eqref{kappa} be fulfilled. The set of such Sonine kernels   is denoted by $\mathcal{L}_{1}$:
\begin{equation}
\label{Son_1}
(\kappa,\, k \in \mathcal{L}_{1} ) \ \Leftrightarrow \ (\kappa,\, k \in C_{-1,0}(0,+\infty))\wedge ((\kappa\, *\, k)(t) \, = \, \{1\}).
\end{equation}
\end{definition}

The properties of the GFI and GFD with the kernels from $\mathcal{L}_{1}$ including the 1st and the 2nd Fundamental Theorems of FC were investigated in \cite{Luc21a}-\cite{Luc21c} on the space $C_{-1}(0,+\infty):=  C_{-1,+\infty}(0,+\infty)$. In particular, we mention the mapping property
\begin{equation}
\label{GFI-map}
\I_{(\kappa)}:\, C_{-1}(0,+\infty)\, \rightarrow C_{-1}(0,+\infty),
\end{equation}
the commutativity law 
\begin{equation}
\label{GFI-com}
\I_{(\kappa_1)}\, \I_{(\kappa_2)} = \I_{(\kappa_2)}\, \I_{(\kappa_1)},\ \kappa_1,\, \kappa_2 \in \mathcal{L}_{1},
\end{equation}
and the index law
\begin{equation}
\label{GFI-index}
\I_{(\kappa_1)}\, \I_{(\kappa_2)} = \I_{(\kappa_1*\kappa_2)},\ \kappa_1,\, \kappa_2 \in \mathcal{L}_{1}
\end{equation}
that are valid on the space $C_{-1}(0,+\infty)$. 

\begin{theorem}[\cite{Luc21a}]
\label{t3}
Let $\kappa \in \mathcal{L}_{1}$ and $k$ be its associated Sonine kernel. 

Then,  the GFD \eqref{FDR-L} is a left inverse operator to the GFI \eqref{FI} on the space $C_{-1}(0,+\infty)$: 
\begin{equation}
\label{FTL}
(\D_{(k)}\, \I_{(\kappa)}\, f) (t) = f(t),\ f\in C_{-1}(0,+\infty),\ t>0,
\end{equation}
and the GFD \eqref{FDC} is a left inverse operator to the GFI \eqref{FI} on the space $C_{-1,(k)}(0,+\infty)$: 
\begin{equation}
\label{FTC}
( _*\D_{(k)}\, \I_{(\kappa)}\, f) (t) = f(t),\ f\in C_{-1,(k)}(0,+\infty),\ t>0,
\end{equation}
where $C_{-1,(k)}(0,+\infty) := \{f:\ f(t)=(\I_{(k)}\, \phi)(t),\ \phi\in C_{-1}(0,+\infty)\}$.
\end{theorem}

\begin{theorem}[\cite{Luc21a}]
\label{t4}
Let $\kappa \in \mathcal{L}_{1}$ and $k$ be its associated Sonine kernel.

Then,  the relations
\begin{equation}
\label{2FTC}
(\I_{(\kappa)}\, _*\D_{(k)}\, f) (t) = f(t)-f(0),\ t>0,
\end{equation}
\begin{equation}
\label{2FTL}
(\I_{(\kappa)}\, \D_{(k)}\, f) (t) = f(t),\ t>0
\end{equation}
hold valid for the functions $f\in C_{-1}^1(0,+\infty) := \{f:\ f^\prime \in C_{-1}(0,+\infty)\}$.
\end{theorem}

In \cite{Luc21a,Luc21b}, the $n$-fold GFIs and GFDs with the Sonine kernels from $\mathcal{L}_{1}$ were defined and studied. For details, we refer the interested readers to these publications. In the rest of this subsection, we discuss the GFI and GFD of arbitrary order that were recently introduced in \cite{Luc21c}. 

To define these operators, we first formulate a condition on their kernels  that generalizes the Sonine condition \eqref{kappa}:
\begin{equation}
\label{Luc}
(\kappa \, * \, k)(t) = \{ 1\}^{<n>}(t),\ n\in \N,\ t>0,
\end{equation}
where
$$
f^{<n>}(t):= (\underbrace{f*\ldots\ * f}_{\mbox{$n$ times}})(t).
$$
Thus, the relation
$$
\{ 1\}^{<n>}(t) = h_{n}(t) = \frac{t^{n-1}}{(n-1)!}
$$
holds valid. Evidently, the Sonine condition corresponds to the case $n=1$ of the more general condition \eqref{Luc}.

Then we proceed with defining the kernels that satisfy the condition \eqref{Luc} and belong to the suitable spaces of functions. 

\begin{definition}[\cite{Luc21c}]
\label{d_c}
Let the functions $\kappa$ and $k$ satisfy the condition \eqref{Luc} and the inclusions $\kappa \in C_{-1}(0,+\infty)$, $k\in C_{-1,0}(0,+\infty)$ hold true.

The set of pairs $(\kappa,\, k)$ of such kernels will be denoted by $\mathcal{L}_n$. 
\end{definition}

As mentioned in \cite{Luc21c}, there are at least two reasonable possibilities to construct a pair $(\kappa_n\, k_n)$ of the kernels from $\mathcal{L}_n,\ n>1$ based on the Sonine kernels $\kappa,\, k$ from $\mathcal{L}_1$. 
The first strategy consists in  building the kernels $\kappa_n = \kappa^{<n>}$ and $k_n=k^{<n>}$. 
The kernels $\kappa_n$ and $k_n$ satisfy  the relation \eqref{Luc} because $\kappa$ and $k$ are the Sonine kernels. However, the pair $(\kappa_n,\ k_n)$ does not always belong to the set $\mathcal{L}_n$. This is the case only  when the inclusion $k^{<n>}\in C_{-1,0}(0,+\infty)$ holds true (of course, $\kappa^{<n>} \in C_{-1}(0,+\infty)$ 
for any $n\in \N$). A more general and important possibility  for construction a pair $(\kappa_n,\, k_n)$ of the kernels from $\mathcal{L}_n,\ n>1$ based on the Sonine kernels $\kappa,\, k$ from $\mathcal{L}_1$ is presented in the following theorem:

\begin{theorem}[\cite{Luc21c}]
\label{tkernel}
Let $(\kappa,\ k)$ be a pair of the Sonine kernels from $\mathcal{L}_1$.  

Then the  pair $(\kappa_n,\ k_n)$ of the kernels given by the formula
\begin{equation}
\label{Son_Luc_n}
\kappa_n(t) = (\{1\}^{<n-1>}\, *\, \kappa)(t),\ \ k_n(t)=k(t)
\end{equation}
belongs to the set $\mathcal{L}_n$.
\end{theorem}

It is worth mentioning that in \cite{Mathematics2021-3}, some important sub-classes of the kernels $\kappa,\, k$ from $\mathcal{L}_n$ in form of convolutions of several different kernels from $\mathcal{L}_1$ were suggested. 

Now we introduce the general fractional integrals and derivatives of arbitrary (non-integer) order and discuss some of their basic properties (for more properties and examples see \cite{Luc21c}).   

\begin{definition}[\cite{Luc21c}]
\label{dao}
Let  $(\kappa,\ k)$ be a pair of the kernels from $\mathcal{L}_n$. 

The GFI with the kernel $\kappa$ has the same form as before
\begin{equation}
\label{GFIn}
(\I_{(\kappa)}\, f)(t) :=  \int_0^t \kappa(t-\tau)f(\tau)\, 
d\tau,\ t>0,
\end{equation}
whereas the GFDs of the Riemann-Liouville and Caputo types with the kernel $k$ are defined as follows:
\begin{equation}
\label{FDR-Ln} 
(\D_{(k)}\, f)(t) := \frac{d^n}{dt^n}\, \int_0^t k(t-\tau)f(\tau)\, d\tau,\ t>0,
\end{equation}
\begin{equation}
\label{FDCn}
( _*\D_{(k)}\, f)(t) := \left(\D_{(k)}\,\left(  f(\cdot) - \sum_{j=0}^{n-1}f^{(j)}(0)h_{j+1}(\cdot)\right)\right)(t),\ t>0. 
\end{equation}
\end{definition}

Evidently, the  GFI \eqref{GFIn} with the kernel $\kappa(t) = h_\alpha(t),\ \alpha>0$ is reduced to the Riemann-Liouville fractional integral \eqref{RLI} and the Riemann-Liouville and Caputo fractional derivatives of the order $\alpha,\ n-1<\alpha<n,\ n \in \N$ are particular cases of the GFDs \eqref{FDR-Ln} and \eqref{FDCn} with the kernel $k(t) = h_{n-\alpha}(t)$. 

To justify the denotation GFI and GFD, we provide formulations of the 1st and 2nd fundamental theorems of FC for the GFDs \eqref{FDR-Ln} and \eqref{FDCn} of the Riemann-Liouville and Caputo types (for the proofs see \cite{Luc21c}).  

\begin{theorem}[\cite{Luc21c}]
\label{t3_n}
Let  $(\kappa,\ k)$ be a pair of the kernels from $\mathcal{L}_n$. 

Then,  the GFD \eqref{FDR-Ln} is a left inverse operator to the GFI \eqref{GFIn} 
on the space $C_{-1}(0,+\infty)$: 
\begin{equation}
\label{FTLn}
(\D_{(k)}\, \I_{(\kappa)}\, f) (t) = f(t),\ f\in C_{-1}(0,+\infty),\ t>0,
\end{equation}
and the GFD \eqref{FDCn} is a left inverse operator to the GFI \eqref{GFIn} on the space 
$C_{-1,(k)}(0,+\infty)$: 
\begin{equation}
\label{FTCn}
( _*\D_{(k)}\, \I_{(\kappa)}\, f) (t) = f(t),\ f\in C_{-1,(k)}(0,+\infty),\ t>0,
\end{equation}
where the space $C_{-1,(k)}(0,+\infty)$ is defined as in Theorem \ref{t3}.
\end{theorem}

\begin{theorem}[\cite{Luc21c}]
\label{t4_n}
Let  $(\kappa,\ k)$ be a pair of the kernels from $\mathcal{L}_n$. 

Then,  the relation
\begin{equation}
\label{2FTCn}
(\I_{(\kappa)}\, _*\D_{(k)}\, f) (t) = f(t) - \sum_{j=0}^{n-1} f^{(j)}(0)\, h_{j+1}(t)
\end{equation}
holds valid on the space $C_{-1}^n(0,+\infty):= \{f:\, f^{(n)} \in C_{-1}(0,+\infty)\}$ and the formula
\begin{equation}
\label{2FTLn}
(\I_{(\kappa)}\, \D_{(k)}\, f) (t) = f(t),\ t>0
\end{equation}
is valid for the functions $f\in C_{-1,(\kappa)}(0,+\infty)$.
\end{theorem}

Further properties and examples of the GFIs and GFDs of arbitrary order were presented in the recent papers  \cite{Luc21c,Mathematics2021-3}. For the results regarding the fractional ODEs and PDEs with the GFDs we refer the readers to \cite{Koch11,Luc21b} and to the recent survey \cite{LucYam20} (see also the references therein). 

\subsection{Open problems}
\label{sec2.5}

In this subsection, we formulate some open problems regarding the one-parameter families of the fractional integrals and derivatives and the GFIs and GFDs.

\vspace{0,2cm}

\begin{enumerate} 

\item 
Deduce a characterization of the fractional integrals for the functions defined on a finite interval that satisfy the desiderata presented in the previous subsection (the Riemann-Liouville fractional integrals are one of such families that satisfies essentially stronger conditions compared to those from desiderata, see Theorem \ref{t_1}).

\item
Characterize all families of the fractional derivatives for the functions defined on a finite interval that are associated to the Riemann-Liouville fractional integrals through the 1st Fundamental Theorem of FC (the $n$th level fractional derivatives including the Riemann-Liouville, Caputo, and Hilfer derivatives are examples of such families).

\item
Deduce a general form of the fractional integrals and derivatives for the functions defined on a finite interval in the sense of the desiderata formulated in the previous section.

\item
Consider the problems formulated above in the case of the functions defined on the positive real semi-axis and on $\mathbb{R}$.

\item
Consider the problems formulated above in the multi-dimensional case, i.e., for the functions defined on $\mathbb{R}^n$.

\item
For a function $\kappa \in C_{-1}(0,+\infty)$, determine the conditions for the inclusion $\kappa \in \mathcal{L}_1$ and derive a representation for the associated kernel $k$.

\item
For a function $\kappa \in C_{-1}(0,+\infty)$, determine the conditions for the inclusion $\kappa \in \mathcal{L}_n$  and derive a representation for the associated kernel $k$.

\item
Adjust the notions of the GFIs and GFDs to the case of a finite interval and $\R$.

\item
Develop a theory of the GFIs and GFDs on other conventional spaces of functions including $L_p$ and H\"older weighted spaces.

\item
Consider the time-fractional differential equations with the GFDs of arbitrary order to model, for instance, the processes intermediate between diffusion and wave propagation.

\end{enumerate}


\section{Numerical Aspects}\label{s3}

For the numerical handling of fractional differential and integral operators, and in particular for the 
numerical solution of initial value problems associated to fractional 
ordinary differential equations, i.e.\ problems of the form 
\begin{equation}
	\label{eq:ivp}
	{}^{\text C} \! D_{a+}^\alpha y(t) = f(t, y(t)), \qquad
	y^{(k)}(a) = y_0^{(k)} \quad (k = 0, 1, \ldots \lceil \alpha \rceil - 1)
\end{equation}
(where the differential operator is chosen to be of Caputo's type \cite[Chapter 3]{Di2010} 
because this is the type that is most important in the mathematical modeling
of processes in physics, engineering, economics, etc.), 
quite a lot of approaches and techniques have
been proposed. Recent surveys about the state of the art in this respect may be found, e.g.,
in \cite{BDST2016,Di2018,Di2019,Ga2018,LZ2015}. In addition to these well known results,
we can witness the publication of a large number of new papers dealing with such aspects.
Unfortunately, many of these papers merely describe the application of well known concepts to
special types of equations, thus not really providing any innovative steps in this area of
research, or---as explained in detail in \cite{DGS2020}---they do not take into account the 
analytic properties of the exact solutions to the differential equations in an appropriate way,
thus stating misleading results.

For the most important numerical approaches, a detailed understanding of their
properties is available (cf., e.g., \cite{DFF2002,DFF2004,Lu1985,Lu1986}).
This includes, in particular, thorough discussions of the weaknesses and limitations
of these algorithms \cite{DFFW2006,FS2001,Ga2010}. In this context, we believe that the
following fundamental issues require attention in the near future and that the corresponding
questions will provide ample opportunities for further research:
\begin{itemize}
\item \emph{The non-locality of fractional differential operators and its implications}: This includes, in
	particular, the observation that traditional numerical algorithms require an $O(N^2)$ amount
	of run time and an $O(N)$ amount of memory when taking $N$ time steps, 
	both of which may be prohibitively large
	in practical application scenarios. Regarding the run time aspect, a number of approaches have
	been developed that provide a remedy. In particular, a general FFT based technique for
	numerically dealing with convolution integrals has been developed in \cite{HLS1985,HLS1988}.
	Its concrete implementation in combination with the traditional Adams method from 
	\cite{DFF2002,DFF2004} has been described in \cite{Ga2018}, thus reducing the 
	computational cost to $O(N \log^2 N)$ (but not lowering the memory requirement). 
	Almost the same can be said about the nested mesh technique of \cite{FS2001} which has an
	$O(N \log N)$ complexity but does not reduce the memory demands either. The approaches of 
	\cite{DF2006,MBSS2015,YY2022} address both matters; they have an $O(N \log N)$ computational 
	complexity and require an $O(\log N)$ amount of memory.
	Here now (to be precise: in Subsections \ref{subs:iss} and \ref{subs:spectral}, respectively), 
	we shall indicate two completely different potential approaches to deal with this
	challenge, neither of which is based on using traditional one-step or
	multistep methods.
	Specifically, the approach of Subsection \ref{subs:iss} has the goal of reducing the 
	complexity even further to the best possible values, i.e.\ $O(N)$ for the computational cost
	and $O(1)$ for the memory, whereas the methods of Subsection \ref{subs:spectral}
	have the goal of providing the approximate solution over the entire interval of interest in one single
	(but large) operation rather than progressing through the interval step by step.
\item \emph{Terminal value problems}: For initial value problems like \eqref{eq:ivp} and their
	numerical handling, a well developed theory is available. However, the application of such mathematical
	models in practice requires the user to be able to provide the initial values $y^{(k)}(a)$, i.e.\ to
	have access to information about the exact solution at time $t = a$, the starting point of the 
	process (indicated by the fact that $a$ is also the starting point of the differential operator 
	${}^{\text C}\! D_{a+}^\alpha$ in the differential equation). In concrete technical applications,
	this information may not necessarily be available; rather, one can sometimes only measure $y(t)$
	(and, if necessary, its derivatives) at some point $t = b > a$. This problem requires numerical
	methods that are completely different from those for \eqref{eq:ivp}; it will be discussed in 
	Subsection \ref{subs:tvp}.
\end{itemize}

\subsection{Algorithms Based on Infinite State Representations}
\label{subs:iss}

In view of the observations recalled in the introduction to 
Section \ref{s3}, a number of novel numerical methods for solving fractional 
differential equations have been developed; see, e.g., 
\cite{Di2008,Di2009,Ba2020,BS2010,Ch2005,HSL2019,KW2021,McL2018,SG2006,SC2006,TR2002,YA2002,ZCSHN2020}.
Although there are various differences in the details, all these methods share the common feature
that they are based on some non-classical representation of the Caputo-type fractional differential operator 
${}^{\text C} \! D_{a+}^\alpha$ of order $\alpha > 0$ that appears in the 
differential equations under consideration, i.e.\ instead of one of the traditional 
forms
\begin{subequations}
\begin{equation}
	\label{eq:def-c1}
	{}^{\text C} \! D_{a+}^\alpha y(t) 
	= \frac 1 {\Gamma(\lceil \alpha \rceil - \alpha)} 
		\int_a^t (t - \tau)^{\lceil \alpha \rceil - \alpha - 1} y^{(\lceil \alpha \rceil)}(\tau) \, \mathrm d \tau
\end{equation}
or
\begin{equation}
	\label{eq:def-c2}
	{}^{\text C} \! D_{a+}^\alpha y(t) 
	= \frac 1 {\Gamma(\lceil \alpha \rceil - \alpha)} 
		\frac{\mathrm d^{\lceil \alpha \rceil}}{\mathrm d t^{\lceil \alpha \rceil}}
		\int_a^t (t - \tau)^{\lceil \alpha \rceil - \alpha - 1} \Bigg[ y(\tau) 
			- \sum_{k=0}^{\lceil \alpha \rceil - 1} \frac{ y^{(k)}(a)}{k!} \tau^k \Bigg] \, \mathrm d \tau ,
\end{equation}
\end{subequations}
one uses a relation of the type \cite[\S 3.2]{Di2010}
\begin{subequations}
\label{eq:rep-nonclass}
\begin{equation}
	\label{eq:int-nonclass}
	{}^{\text C} \! D_{a+}^\alpha y(t) 
	= \int_0^\infty \phi(w, t) \, \mathrm d w,
\end{equation}
the so-called infinite state representation (or diffusive representation, cf.\ \cite{Mo1998}) 
of the fractional derivative of $y$. Here the integrand $\phi$ 
(whose values $\phi(w, t)$ for $w \in (0, \infty)$ are known as the
infinite states of the observed system at time $t$) solves an inhomogeneous linear first order 
initial value problem of the form
\begin{equation}
	\label{eq:ivp-nonclass}
	\frac{\partial}{\partial t} \phi (w, t) 
	= h_1(w) \phi(w,t) + h_2(w) y^{(\lceil \alpha \rceil)}(t), 
	\qquad
	\phi(w, 0) = 0,
\end{equation}
with certain functions $h_1, h_2 : (0, \infty) \to \mathbb R$ for which many different specific choices
have been suggested in the papers mentioned above.
\end{subequations}
Using this representation, the task of numerically evaluating the fractional differential operator
(which is the core of every numerical method for solving such differential equations) amounts 
to
\begin{enumerate}
\item approximately computing the solution $\phi(w, t)$ to the initial value 
	problem \eqref{eq:ivp-nonclass} at the
	time point $t = t_n$ currently under consideration for certain suitably chosen values of $w$, and
\item based on this information, 
	numerically evaluating the integral on the right-hand side of eq.\ \eqref{eq:int-nonclass}
	with a properly designed quadrature formula.
\end{enumerate}
In this context, the ``suitable choice'' of the values $w$ mentioned in the first of the two items
above is essentially determined by the location of the nodes of the quadrature formula used in the
second item. 

From an algorithmic and software engineering point of view, this is an excellent approach because
it is both fast---having a computational complexity of $O(N)$ when $N$ time steps need to be taken
rather than the $O(N^2)$ of the original fractional Adams method from \cite{DFF2002} or the
at best $O(N \log N)$ cost of its modified version mentioned above---and memory efficient
in the sense that a computation over $N$ time steps requires only a memory amount of $O(1)$ 
instead of $O(\log N)$ as in the method from \cite{DF2006,MBSS2015,YY2022} or even $O(N)$
in all other versions of the Adams method (and, in this comparison, replacing the Adams method 
by a different classical approach would not make any difference in these respects \cite{Di2019}).
This observation is an immediate consequence of the fact that one has to solve the differential
equation \eqref{eq:ivp-nonclass} which is a differential equation of order 1 and hence
does not exhibit any non-locality or any other memory effects.

However, from a numerical analysis perspective, the situation is highly unsatisfactory for a number
of reasons:
\begin{itemize}
\item Almost all the arguments on which the specific versions of the algorithm indicated above are
	based are heuristic in nature. In particular, the quadrature formulas for the integral from 
	eq.~\eqref{eq:int-nonclass} required in step~2 are constructed in an ad hoc manner; the analytical
	properties of the integrand that are well known to strongly influence the accuracy
	of the approximate results \cite{BP2011} are taken into account only partially in the 
	proposal from \cite{Di2009} and are completely ignored in most of the other schemes.
\item There is hardly any investigation regarding the interaction between the numerical ODE solver
	required in step~1 of the algorithm and the quadrature formula from step~2. The only 
	information usually provided is the observation that the nodes of the quadrature
	formulas that have been proposed have such magnitudes that the solution of the 
	ODE requires an A-stable method. Then, according to Dahlquist's well known second
	barrier \cite[\S V.1]{HW2002}, the order of the ODE solver cannot be higher than 2 if it is chosen
	from the set of linear multistep methods; a higher order can however be achieved by using
	suitable implicit Runge-Kutta methods \cite[\S IV.6]{HW2002}. In practice, to avoid the intricacies
	connected to higher order Runge-Kutta methods, all the proposals found in the literature suggest
	to use a linear multistep method (often, the backward Euler method). 
	The low order seems to be acceptable because, as mentioned above, no strong effort has 
	been invested into the
	selection of the quadrature scheme, so using a high order but complex ODE solver
	does not seem to be worth while anyway.
\item The typically chosen building blocks for the overall algorithm contain a number of parameters
	(number and location of quadrature nodes, quadrature weights, choice of
	the ODE solver, step size of ODE solver) whose influence on the accuracy of the final result
	is essentially unknown. None of the papers mentioned above contains any comprehensive error
	analysis.
\item For many of the approaches mentioned above, it is not immediately clear how the 
	implied constants in the error bounds (if such bounds are known at all) depend on the 
	order $\alpha$ of the differential operators in question; in particular, it is often unknown 
	whether the bounds hold uniformly for all $\alpha \in (0,1)$ or whether the constants
	blow up as $\alpha \to 0$ or $\alpha \to 1$.
\end{itemize}
Therefore, there exists a strong demand for a thorough investigation of the non-classical algorithms
for the numerical solution of fractional differential equations.

\subsection{Spectral Methods and Related Techniques}
\label{subs:spectral}

A frequent objection against the use of one-step or multistep methods for solving fractional
differential equations is that these methods are of some kind of a local nature which does not
appear to be a natural concept for dealing with operators like fractional derivatives that are
not local. One potential and promising approach in this direction is the use of spectral methods.

Traditionally, such spectral methods are based on classical polynomials \cite{LX2009,LZL2012},
but clearly this is not ideal because polynomials cannot accurately capture the singular
behaviour that the exact solutions to fractional differential equation exhibit near the 
starting point of the associated differential operator \cite{SS2019}. A more suitable 
alternative is to use generalized polynomials with a combination of integer and 
non-integer exponents as suggested, 
e.g., in \cite{GSW2009,ZK2013,ZK2014}. This alternative comes with two significant
advantages: 
\begin{enumerate}
\item The exponents can be chosen to exactly match those present in the analytical
	solution to the differential equation, thus eliminating the error component 
	induced by these (generalized) monomials and leading to a significantly better
	convergence rate.
\item If the basis of the set of generalized polynomials is chosen in an appropriate
	way \cite{ZK2013} then, for certain classes of differential equations, the coefficient
	matrices of the discretized differential equations become sparse, thus greatly
	reducing the computational cost and the memory requirements of the 
	resulting algorithms.
\end{enumerate}

One of the key outcomes of \cite{ZK2014} is the recommendation to use a 
Petrov-Galerkin strategy for the construction of the approximate solution, i.e.\ one
should choose the test function space for the Galerkin method to be different
from the space of trial functions in which the approximate solution is sought.
Following the development of \cite{ZK2014}, we shall summarize this
construction by means of the simplest possible example application, 
namely the initial value problem
\begin{equation}
	\label{eq:ivp-spectral}
	 {}^{\mathrm C} \! D_{-1+}^\alpha y(t) = f(t), \qquad y(-1) = y_0, \quad t \in [-1,1]
\end{equation}
i.e.\ a differential equation whose right-hand side is independent of the 
unknown function $y$ so that the analytical solution process essentially
amounts to a fractional integration of the forcing function.
(The interval $[-1,1]$ has been chosen for convenience and without loss of generality.)
In this special case, one starts with the construction of an associated 
Sturm-Liouville eigenvalue problem of the first kind \cite{ZK2013},
\begin{eqnarray}
	\nonumber
	{}^{\mathrm{RL}} \! D_{1-}^{\alpha/2} K \, {}^{\mathrm C} \! D_{-1+}^{\alpha/2} \Phi_{1, \lambda} (t)
	+ \lambda (1-t)^{-\alpha/2} (1+t)^{-\alpha/2} \Phi_{1, \lambda} (t) &=& 0, \\
	\label{eq:slevp1}
	\Phi_{1, \lambda}(-1) &=& 0, \\
	\nonumber
	{}^{\mathrm{RL}} \! I_{1-}^{1-\alpha/2} K \, 
		{}^{\mathrm C} \! D_{-1+}^{\alpha/2} \Phi_{1, \lambda} (t) \Big|_{t = 1} &=& 0
\end{eqnarray}
with some constant $K$, and computes the eigenfunctions of this problem. It turns out that
these eigenfunctions can be written as
\begin{equation}
	\label{eq:eigenf1}
	\mathcal P_{1,n}(t) = (1+t)^{\alpha/2} P_{n-1}^{-\alpha/2, \alpha/2}(t)
	\qquad (n = 1, 2, 3, \ldots)
\end{equation}
where $P_k^{(-\alpha/2, \alpha/2)}$ denotes the standard $k$-th degree Jacobi polynomial 
for the weight function defined by the superscript in the usual way \cite[Chapter IV]{Sz1975}. 
The eigenfunctions
$\mathcal P_{1,n}$ for $n \in \{1, 2, \ldots, N\}$ with some predefined natural number $N$
are chosen as the basis functions of the space in which the approximate solution 
is sought.

Moreover, a similar Sturm-Liouville eigenvalue problem of the second kind \cite{ZK2013}
\begin{eqnarray}
	\nonumber
	{}^{\mathrm{RL}} \! D_{-1+}^{\alpha/2} K \, {}^{\mathrm C} \! D_{1-}^{\alpha/2} \Phi_{2, \lambda} (t)
	+ \lambda (1-t)^{-\alpha/2} (1+t)^{-\alpha/2} \Phi_{2, \lambda} (t) &=& 0, \\
	\label{eq:slevp2}
	\Phi_{2, \lambda}(1) &=& 0, \\
	\nonumber
	{}^{\mathrm{RL}} \! I_{-1+}^{1-\alpha/2} K \, 
		{}^{\mathrm C} \! D_{1-}^{\alpha/2} \Phi_{2, \lambda} (t) \Big|_{t = -1} &=& 0
\end{eqnarray}
is constructed, and the associated eigenfunctions of this problem are determined too.
In this case, they can be written in the form
\begin{equation}
	\label{eq:eigenf2}
	\mathcal P_{2,n}(t) = (1-t)^{\alpha/2} P_{n-1}^{\alpha/2, -\alpha/2}(t)
	\qquad (n = 1, 2, 3, \ldots),
\end{equation}
and these functions (for $n = 1, 2,\ldots, N$) are used as the test functions.

Based on these definitions of trial and test functions, one can then apply the
standard Petrov-Galerkin scheme to compute the approximate solution $Y_{N,0}$ of the 
differential equation from the given problem \eqref{eq:ivp-spectral} augmented
with the initial condition $y(-1) = 0$. The fact that we can only work with this initial
condition is due to the observation that, as can be seen from eq.~\eqref{eq:eigenf1},
all basis functions satisfy $\mathcal P_{1,n}(-1) = 0$, and so the approximate 
solution---being a linear combination of these basis functions---must also satisfy
$Y_{N,0}(-1) = 0$. The case of the general initial condition given in eq.~\eqref{eq:ivp-spectral}
then leads to the approximate solution $Y_N(t) = y_0 + Y_{N,0}(t)$ with the function
$Y_{n,0}$ from above being the solution of the differential equation with a homogeneous 
initial condition. This follows from the fact that the additive constant $y_0$ is annihilated
by the Caputo differential operator in the differential equation \eqref{eq:ivp-spectral}.

A particularly relavant feature of this process is that the stiffness matrix of the resulting linear system
that determines the coefficients of the approximate solution $Y_{N,0}$ is actually a diagonal
matrix, so the system can be solved in an extremely simple way.

From a more general perspective, one starts from the given differential 
equation and constructs associated Sturm-Liouville eigenvalue problems of the
first kind and of the second kind, respectively. The eigenfunctions 
of the eigenvalue problem of the first kind are then used the basis functions
for the space of approximate solutions, and the eigenfunctions of the
eigenvalue problem of the second kind are chosen as the test functions of the
Petrov-Galerkin algorithm.

A detailed and recent survey describing the state of the art in this area of research
may be found in \cite{SS2019}.

In a similar way, one can handle so-called spectral element methods where the
domain is first decomposed in a finite element like manner, and afterwards
a spectral kind of approximation is used on each element; cf., e.g., the
survey in \cite{LZZ2019}.

The description given above indicates that the strategy proposed in \cite{ZK2014}
depends on the precise form of the given differential equation in various respects.
This applies, e.g., to the construction of the two Sturm-Liouville eigenvalue problems
and the determination of their eigenfunctions (i.e.\ the trial and test functions), 
but also to the specific way in which a
potential inhomogeneous initial condition can be incorporated (although the latter
can be generalized rather easily). 
From a practical perspective, it would be advantageous to have a set of test functions
and a set of trial functions that could be used not just for one type of differential
equations but for a broad class of such problems. The search for sets of such 
basis functions (that, of course, should be able to guarantee a rapid convergence
of the sequence of approximate solutions towards the exact solution) is a topic
of great interest. Moreover, one might raise the question whether a different construction
of the basis and test functions could lead to an even faster convergence.

A notable weakness of this approach that might also be addressed in future research
projects is that one is naturally lead to consider the
covergence behaviour and the error estimates in an $L_2$ norm only whereas other
solvers for fractional ordinary differential equations are usually analyzed with 
respect to the stronger $L_\infty$ norm.

\subsection{Terminal Value Problems}
\label{subs:tvp}

A completely different but also potentially very relevant question with respect to
which many open questions exist is also connected to the fractional differential equation 
\begin{subequations}
\label{eq:tvp}
\begin{equation}
	\label{eq:fde1}
	{}^{\text C} \! D_{a+}^\alpha y(t) = f(t, y(t))
\end{equation}
that is discussed with some given function $f$ on the interval $[a, a + T]$ for some $T>0$.
We concentrate on the practically most important case $\alpha \in (0,1)$. Keeping in mind the 
well kown properties of the operator ${}^{\text C} \! D_{a+}^\alpha$ and, in particular,
the classical existence and uniqueness result for initial value problems \cite[Theorem 6.5]{Di2010},
it is natural to believe that this problem requires exactly one additional condition of the form 
\begin{equation}
	\label{eq:tc1}
	y(b) = y^*
\end{equation}
\end{subequations}
with some given numbers $y_0 \in \mathbb R$ and $b \in [a, a+T]$
to assert the uniqueness of the solution. Such a condition with $b > a$ may be
necessary for example because
measurements of the state of the modeled system are available only at time $t = b$ but not
at time $t = a$ (the starting point of the differential operator in eq.~\eqref{eq:fde1}).
It turns out in the case $b > a$, however, that the uniqueness can only be guaranteed in the
case where the problem is scalar, i.e.\ when we are looking for a function $y : [a, a+T] \to \mathbb R$,
but not in the vector-valued case where $y : [a, a+T] \to \mathbb R^d$ with $d \ge 2$, cf.\ \cite{CT2017}.
We thus impose the additional restriction to consider only the scalar case. In this situation,
the problem can be decomposed into two steps:
\begin{enumerate}
\item Solve the problem on the interval $[a, b]$.
\item Solve the problem on the interval $[a, a + T]$.
\end{enumerate}
Once the first step has been completed, the initial value $y(a)$ is known, and hence
it can be added to the differential equation in the form of an initial condition, replacing
the terminal condition~\eqref{eq:tc1} and thus creating an initial 
value problem whose solution can be computed with standard and well established methods to 
find the solution of the second step and hence the overall solution. Therefore, the only really
novel task is step 1. Since the auxiliary condition is formulated at the end point of the
interval of interest, this problem is called a \emph{terminal value problem}.

From an analytical point of view, the fundamental difference between initial and terminal
value problems is highlighted when the two problems are rewritten as equivalent integral
equations: While, as is well known, an initial value problem is equivalent to a 
\emph{Volterra} integral equation \cite[Lemma 6.2]{Di2010}, the terminal value 
problem \eqref{eq:tvp} can be reformulated on $[a,b]$ as
\begin{subequations}
\begin{equation}
	\label{eq:fredholm}
	  y(t) =  y^* + \frac 1 {\Gamma(\alpha)} \int_a^b G(t,s) f(s, y(s)) \, {\mathop{\text d}} s 
\end{equation}
with
\begin{equation}
	\label{eq:green}
	G(t, s)  
    		= \begin{cases}
		      - (b - s)^{\alpha-1} &  \mbox{ for } s > t, \\
		        (t - s)^{\alpha-1} - (b - s)^{\alpha-1} & \mbox{ for } s \le t,
		    \end{cases}
\end{equation}
\end{subequations}
which is a \emph{Fredholm} integral equation \cite[Theorem 6.18]{Di2010}. From this 
difference in the types of the integral equations, it is already clear that the terminal value 
problem needs completely different numerical methods than the initial value problem; 
in particular, one can see that the nature of a fractional order terminal value problem
is much closer to an integer order boundary value problem than to an integer order 
initial value problem. Therefore, it is natural to solve fractional order terminal value problems
with the help of methods based on the construction principles coming from the area
of integer order boundary value problems. 

Following the line of thought outlined above, a first obvious strategy to solve terminal
value problems is the class of shooting methods (cf., e.g., \cite[Chapter 2]{Ke1992} for a
very detailed discussion of such methods in the case of integer order equations): One starts 
with an initial guess for the initial value $y(a)$, solves the associated initial value problem, 
and compares the terminal value, i.e.\ the approximation for $y(b)$, obtained in this way with
the desired value $y^*$. If these two values are sufficiently close to each other, the approximate
solution is accepted; otherwise, a new attempt is started with a suitably modified initial value. 

Such methods have been proposed for fractional terminal value problems in \cite[\S 6]{DF2012} 
where the classical Adams method from \cite{DFF2002,DFF2004} has been suggested for 
solving the initial value problems. Various improvements of this algorithm with the aim of increasing
its efficiency have been discussed in \cite{Di2015}. One of the main observations from that paper is
that in the early stages of the iterative process, i.e.\ when the selected initial value is still far away
from the correct one, it is not necessary to use a very small step size for the initial value problem
solver. However, that proposal still leaves room for improvement in this direction, and numerous 
other questions remain open as well, most importantly the question what the best algorithm for
selecting the initial value in the next iteration step is. This question is currently under active 
consideration \cite{DU2021}.

A case study conducted by Ford and Morgado \cite{FM2011} provides a comparison of this
shooting method based on the Adams scheme with corresponding approaches using other
solvers for the initial value problems, essentially based on
linear multistep methods, and comes to the conclusion that, among the methods
investigated, the one using Adams algorithm provides the best balance between accuracy
and computational effort. Obviously, this investigation can only cover a limited
number of alternatives to the Adams method and might be extended to other schemes.
A different replacement for the Adams scheme, namely a collocation method using 
a combination of classical piecewise polynomials and generalized polynomials (with 
non-integer exponents) has been shown to also provide relatively good results in \cite{FMR2014}.

Shooting methods are, of course, not the only possible approach to numerically handle 
fractional terminal value problems. In view of the integral equation represtentation 
\eqref{eq:fredholm} of the problem, all the known standard techniques for the numerical
solution of weakly singular Fredholm equations are available too, cf., e.g., \cite{Ha1995}.
It seems, however, that among the many classes of methods that can be used in principle
such as, e.g., Nystr\"om methods, collocation methods, or Galerkin methods, only a very
small number has actually been attempted for fractional terminal value problems. A notable
instance is a collocation scheme, again based on a combination of classical polynomials 
and polynomials with non-integer exponents (which, in view of the well known 
properties of the exact solutions to fractional differential equations, seems to be a 
quite natural idea) discussed in \cite{FMR2015}. The convergence results shown in that 
paper indicate the method to be quite attractive and rapidly convergent. A comparison of the 
computational cost of this method to the cost of a shooting method would be useful to 
get a comprehensive impression. Similar studies for other approaches to the Fredholm
equation (Galerkin, Nystr\"om, etc.) would be most welcome. Such investigations should
in particular include a discussion about a suitable set of basis functions that combines
the ability to accurately model the asymptotic behaviour of the exact solution as $t \to a$,
thus exactly capturing the singular behaviour of its derivatives, with a rapid convergence
of the remaining part of the approximate solution to the smooth part of the exact solution.
This latter aspect in particular includes the task of applying the spectral or spectral element
methods discussed in Subsection \ref{subs:spectral} to the general class of terminal value
problems described in the present subsection.



\section{Applications in Physics}\label{s4}



For the first time, FC was applied by Niels Henrik Abel in 1823 \cite{Abel1823,Abel2012} to the tautochrone problem, which is the kinematic problem of finding the trajectory of point mass sliding without friction to its lowest point, when time does not depend on its starting point of the trajectory.
Then applications of FC to problems of physics were proposed in the works of Joseph Liouville, beginning in 1832 \cite{Liouville1832,Lutzen1985,Lutzen1982,Lutzen1990}. 
For Liouville, the applications of FC were not easy illustration to awaken the interest of other mathematicians. 
The applicability of FC was the meaning of the existence of this mathematical theory. 
This relationship to FC significantly distinguishes Liouville from his predecessors, 
whose main goal was a purely mathematical generalization of ordinary calculus \cite[p. 307]{Lutzen1990}. 
Liouville's contribution to the application of FC is described by J. L\"utzen in 
\cite{Lutzen1985,Lutzen1982} and \cite[pp. 307--320]{Lutzen1990}. 
The contribution of the first pioneers of the application of FC in physics is described in \cite{Valerio:14}. 


Fractional calculus, equations with fractional integrals, derivatives and differences are a powerful tool 
for describing local processes in time and space with different types of nonlocality. 
Conventionally, all applications can be divided into five directions: 
nonlocality in space and time in the framework of continuous and discrete approaches, and their interrelationships, which can be called 
(CT) Nonlocality in continuous-time;
(DT) Nonlocality in discrete-time; 
(CS) Nonlocality in continuous-space; 
(DS) Nonlocality in discrete-time; 
(R) Relationship between these directions.

Moreover, there is a need to classify the types of nonlocality and to link these types with the description of the types of phenomena. 
For example, it is necessary to distinguish nonlocalities in time that describe (1) memory, (2) distributed lag (time delay); (3) distributed scaling (dilatation, dilation), and others. 
Note that these different phenomena and types of nonlocality must be described by different types of kernels of fractional operators.

One of the most important trends in the modern stage of applications of FC is the urgent need to describe general types of nonlocal phenomena and the corresponding sets of kernels of fractional operators that form FC. We emphasize here the importance of not only, and not so much specific examples from various fields of science, since there can be an infinite number of such examples. It is important to solve general questions, it is important to answer the questions of what types of operator kernels (and therefore fractional operators) what types of phenomena can describe \cite{Mathematics2020}. For example, (1) the Riemann-Liouville and Caputo fractional operators can describe memory; (2) the fractional operators of Hadamard type and 
Erdelyi-Kober can describe scaling and dilation; (3) Fractional operators in the form of a Laplace convolution, whose kernel satisfies the normalization condition, can describe a distributed lag (time delay), and so on \cite{Mathematics2020}.
 
The very important trend in application of FC is investigations and results concerning of general form of nonlocality, which can be described by general form of operator kernels, and not its particular implementations and representations \cite{Mathematics2021-2}. 

Due to page limits for review, we will highlight only applications of FC that seem very promising, important, and interesting. 
Unfortunately, some of the important and promising applications of FC will remain out of consideration. 
Some of these areas are described in volumes 4 and 5 of Handbook of Fractional Calculus with Applications \cite{Tarasov:19a,Tarasov:19b}, which contain 25 reviews on various areas of applications of FC in physics. 
Other areas of application of FC are described in reviews \cite{IJMPB2013,Ionescu,BOOK-MDPI-2020,BOOK-DG-2021} and in books on fractional dynamics \cite{Zaslavsky2004,Springer2010,KlafterLimMetzler2012}.


\subsection{Nonlocal Continuum and Lattice Mechanics}


The classical area of application of FC is the continuum mechanics, which remains in the trend of future research. 
At the same time, the interest is currently shifting from nonlocalities in time to nonlocalities in space, and in the direction of nonequilibrium thermodynamics, and electrodynamics of media with space dispersion and frequency dispersion. 

Continuum mechanics with nonlocality in time has a long history dating back to the work of Ludwig Boltzmann in 1874.
The first physical model of media with nonlocaliy in time has been proposed by Ludwig Boltzmann in 1874 and 1876 \cite{Boltzmann1,Boltzmann2} for isotropic viscoelastic media. In these works, te fading and superposition principles were suggested. 
Then this approach is used in the works of Vito Volterra \cite{Volterra}.
The first application of FC to continuum mechanics of viscoelastic as media with nonlocality in time was proposed in work of Andrey N. Gerasimov \cite{Gerasimov} in 1948. 
At present time, FC is actively used to describe continuum mechanics with nonlocality in time (for example, see \cite{Rabotnov2,Mainardi,APSS2014a,Povstenko2015} and references therein).

Continuous mechanics with nonlocality in space, in contrast to mechanics with nonlocality of time, has a shorter history. 
Apparently the first application of FC to continuum mechanics with nonlocality in space was proposed by V.S. Gubenko in 1957 \cite{Gubenko,Rostovtsev}. The theory of nonlocal continuum mechanics was formally initiated by the papers of Kr\"oner \cite{Kroner}, and Eringen \cite{Eringen1972a,EringenEdelen1972,Eringen1972b}.
Nonlocal theory is based on the assumption that the forces between material points are a long-range type, thus reflecting the long-range character of inter-atomic forces \cite{Eringen2002,Rogula}. 

In electrodynamics and mechanics of media the nonlocality in space is interpreted as special form of spatial dispersion (SD) \cite{SR1961,ABR1978,Agr1984}.
For the first time, FC was used to describe SD in electrodynamics and mechanics in works 
\cite{AP2013,MPLB2016} and  \cite{CEJP2013,ISRN-CMP2014,MOM2014,IJSS2014},  where fractional form of SD are described. 

For continuous models of media with nonlocality in space, the 
limitation is the absence of a fractional analogue of usual tensor and vector calculus \cite{AP2008}, \cite[pp. 241--264]{Springer2010}. 
This is due to the fact that nonlinear coordinate transformations transform fractional operators of non-integer order into pseudo-differential operators, accompanied by a change in the form of nonlocality.
For the first time, a consistent mathematical formulation of fractional vector calculus was proposed in work \cite{AP2008}, 
which includes fractional generalizations of the differential operations (gradient, divergence, curl), the integral operations (flux, circulation), and the relationship of these operators in the form the generalized Gauss, Stokes and Green theorems (see also \cite[pp. 241--264]{Springer2010}).


It is well known in solid state physics that long-range interactions between particles can be considered as a source of nonlocality in space, which largely determines the properties of the medium \cite{Eringen2002}. 
Because of this, it seems natural to construct mechanics of continuum and lattice with nonlocality in space \cite{Eringen2002,Kunin1982,Kunin1983}. 
Long-range interactions in lattice models of media is important in physics, starting with Dyson's work in 1969 \cite{Dyson1,Dyson2} about models of spin chains and lattices with long-range interactions.
However, the active use of FC in this direction begins only now. 

The use of FC for constructing lattice models with nonlocality in space (with long-range interactions) and consideration of the continuous (continuum) limit of these models, which is described by fractional differential equations, actually began with the works \cite{JPA2006,JMP2006,CNSNS2006,Chaos2006} in 2006 (see also \cite[pp. 153--214]{Springer2010}).
Then the proposed approach was generalized from one-dimensional case for $n$-dimensional case, and was represented as a lattice FC
\cite{JPA2014,AMC2015,FCAA2016,AHEP2014,MPLA2021}.
This representation allows to formulate concept of exact differences of integer and fractional orders in \cite{JM2015,PLA2016,M2016} and \cite{CNSNS2016-2,CMA2017}. 
Unfortunately, for exact fractional differences, which are proposed in \cite{CNSNS2016-2,CMA2017,FCAA2016}, the fundamental theorems of FC have not been proven for non-integer orders of operators at present time.


\subsection{Economics with Nonlocality in Time: Memory and Distributed Lag}


The first application of fractional calculus in economics can be considered the work of the Nobel laureate C.W.J. Granger, R. Joyeux, in 1980 \cite{Granger4}, in which the Gr\"{u}nwald-Letnikov fractional differences were actually rediscovered. 
In this work, a fractional generalization of ARIMA models was proposed. 
In recent years, various attempts to apply FC in economics and finance are suggested. 
The history of these applications of FC is given in \cite{History} and \cite[pp.~5--32]{BOOK-MDPI-2020}, \cite{BOOK-DG-2021}.
Many of these attempts have various disadvantages that described in the review \cite{Rules}. 

FC is used to construct the basic concepts of economic theory and to form mathematical models of economic processes on their basis \cite{BOOK-DG-2021}. 
The directions of further applications of FC in this area are described in \cite{History} (or \cite[pp.~5--32]{BOOK-MDPI-2020}), and in Section 30 \cite[pp. 532--539]{BOOK-DG-2021}.
Additional references on application of FC in economics and finance can be found in the list of references in \cite{BOOK-DG-2021,History}. 

Note that, in contrast to mechanics and physics, where relaxation equations are often used, in economics, growth equations are more often applied. The fractional differential equations describing relaxation and growth have a qualitative difference in the behavior of the solution due to their difference in asymptotic behavior. This is one of the distinctive features of the novelty of the behavior of economic processes with memory from processes with memory in physics. 
The appearance of memory or memory changes can lead to a significant acceleration or slowdown of the process, that is, the characteristic time of the process can be substantially increased or decreased by memory. 


Application of FC is important to describe nonlocality in the form of distributed lag (time delay).
The distributed lag has been considered starting with the models: (1) the model with uniform distributed lag that was suggested by Michal A. Kalecki \cite{Kalecki} in 1935 to described business cycles \cite[pp.~251--254]{Allen63}; (2) the models with the exponential distribution of delay time was proposed by Alban W.H. Phillips \cite{Phillips1} in 1954; (3) models and operators with distributed lag that were considered by Roy G.D. Allen \cite{Allen63} in 1956. 
The distributed lag is described by the kernels that are probability density functions (p.d.f), when the operators are defined by the Laplace convolution with these kernels. 
The distributed lag is caused by finite speeds of processes, and cannot interpreted as a memory. In physics, the distributed lag as a form of nonlocality in time is interpreted as hysteresis. 

Simultaneous presence of memory effects and distributed lag is important for economic models. 
However, for fractional operators with distributed lag, which are proposed in \cite{CNSNS-LAG} and used in \cite{Lag-Keynesian1,Lag-Phillips}, \cite[Ch.20-23]{BOOK-DG-2021}, the fundamental theorems of FC have not been proved in general. Note that p.d.f of the gamma distribution is the Sonin kernel, and we can use the general FC \cite{Kochubei2011,Samko2003,Luc21a,Luc21b,Luc21c,Mathematics2021-3} in this case and some other distributions on the positive semiaxis. 
This is an important direction of future research in FC and its applications in economics and physics.


\subsection{Non-Markovian Dynamics of Open Quantum Systems}


The theory of open quantum systems is the most general type of modern quantum mechanics as fundamental theory \cite{Op3,TBook2008}.
Moreover, this theory has great practical importance for the creation of quantum computers and quantum informatics. The influence of the environment changes the quantum computation, which is realired by dynamics of quantum systems of qubits.
Note that the Schr\"{o}dinger equation describes only pure quantum states and cannot be used for general description of quantum dynamics. 

In recent decades, the theory of open quantum systems has been actively developing (for example, see basic papers \cite{K1,Lind1,Lind2}, books \cite{Op3,Op5,TBook2008}, and reviews \cite{ISSSS}). 
The most general form of Markovian equations, which describe quantum observables and quantum states of quantum systems, was suggested by Gorini, Kossakowski, Sudarshan and 
Lindblad in \cite{K1,Lind1,Lind2}. 
The most general explicit form of equations for the quantum mechanical systems is the Lindblad equations, which describe Markovian dynamics of quantum observables and quantum states.

We should emphasize that the openness of classical and quantum systems (the interaction with the environment) can lead to nonlocality in time \cite{CEJP2012}. 

Currently, modern quantum mechanics is faced with the question of the most general form of the equation describing the non-Markovian dynamics of quantum systems with memory. 
Attempts to construct a non-Markovian theory of open quantum systems with memory have been actively pursued recently (for example, \cite{QUANT-MEM-1,QUANT-MEM-2,QUANT-MEM-3,QUANT-MEM-4,QUANT-MEM-5,QUANT-MEM-6}, and references therein).
The non-Markovian character of quantum processes was often interpreted as nonlocality in time. 
Note that exact solutions of equations describing non-Markovian dynamics with memory have practically not been proposed in these works. 
All these attempts were not associated with the use of FC. 
In construction of a consistent non-Markvian theory of open quantum systems, the FC gives new possibilities for the development of this direction. 

For the first time, the use of FC to take into account memory effects (non-Markovity) in open quantum systems was proposed in work \cite{Springer2010} (see Chapter 20 in book 
\cite[pp. 477--482]{Springer2010} and \cite{AP2012,Handbook-5-Open,Entropy2021}).
The solutions of non-Markovian equations describing quantum systems with memory, which is based on FC, were proposed in these works. 

For the first time, fractional powers of Lindblad superoperators were defined and used to describe non-Markovian dynamics in 2008 \cite{TBook2008} (see Chapter 20 in book \cite[pp. 433--444]{TBook2008}, Chapter 20 in book \cite[pp. 458--464, 468--477]{Springer2010}, \cite{TMP2009,FD-RA2012}). 
Solutions of the proposed generaized Lindblad equations, which describe non-Markovian quantum dynamics, were derived in \cite{TBook2008,TMP2009,FD-RA2012,Springer2010}.

For construction of a theory of non-Markovian quantum dynamics, it is important to consider time-dependent parameters \cite{AP2017,Handbook-5-Open} and the relativistic quantum systems with memory \cite{PLA2020}.
It is important for the non-Markovian quantum dynamics the generalization in the directions: the path integral for open quantum systems \cite{JPA2004}, \cite[pp. 475--485]{TBook2008};
uncertainty relation for open quantum systems \cite{JMP2013}; the pure stationary states of open quantum systems \cite{PRE2002,PLA2002}, \cite[pp. 453--462]{TBook2008}.
The non-Markovian quantum dynamics can also be considered in the framework of fractional generalizations of Heisenberg equations \cite{PLA2008,TBook2008}, \cite[pp. 457--466]{Springer2010}, in contrast to the time-fractional Schr\"{o}dinger equation. 
For further reference, it can be used the list of references in \cite{Entropy2021}. 

Quantum mechanics with nonlocality in space and quantum field theory with non-locality in space-time by using models of lattices with long-range interactions
and its continuum limits \cite{JPA2014,AMC2015,FCAA2016,AHEP2014,MPLA2021,JM2015,PLA2016,M2016,CNSNS2016-2,CMA2017} are important future direction of applications of FC to quantum physics.
This is due to the fact that at present various nonlocal models are being actively investigated within the framework of quantum field theory at present time.
In addition, lattice models with long-range action are being actively investigated within the framework of quantum physics of solids and condensed matter.


\subsection{Nonlocal Discrete Maps} 


In nonlinear dynamics, maps with discrete time are derived from integer-order differential equations with periodic kicks (Sec. 5.2, 5.3 in \cite[pp. 60--68]{Zaslavsky2004}, and Ch. 18 in \cite[pp. 409--453]{Springer2010}). 
These maps define the next step only by the previous step (or a fixed number of previous steps. 
In discrete maps with nonlocality in time means that the next step depends on all past steps.

Discrete maps from fractional differential and integral equations of non-integer orders can be also derived.
For the first time such discrete maps were obtained in article \cite{Tarasov-Zaslavsky,Tarasov-Map1,Tarasov-Map2} in 2008. 
The proposed nonlocal maps can be derived by using the equivalence of the FDE and the Volterra integral equations in \cite{Tarasov-Map1,Tarasov-Map2,Springer2010}.
Then, this approach has been applied in \cite{Tarasov-Map3,TT-Logistic,Entropy2021,MMAS2021,BOOK-DG-2021,Mathematics2021-2,CNSNS2021,CSF2021}. 
The first computer simulations of the suggested nonlocal maps were made in \cite{TT-Edelman1,TT-Edelman2}, and then numerical simulations have proved the existence of new types of attractors and new types of chaotic behavior for these maps (see \cite{Edelman1,Edelman3,Edelman-Handbook-2,Edelman-Handbook-4,Edelman2021} and references therein).

A new general approach to the study of the chaotic behavior of discrete maps with nonlocality was proposed in \cite{Mathematics2021-2,Edelman2021}. In the framework of this approach, new types of attractors and chaotic nonlocal dynamics
can be discovered for nonlocal discrete maps with discrete convolution \cite{Mathematics2021-2,Edelman2021}.


\subsection{Self-Organization with Memory}


Self-organization is a process of formation of ordered spatial or temporal structures that can be realized without external influences. Self-organization can be realized in physical, chemical, biological and economic processes \cite{Self-Org1,Self-Org2,Self-Org3,Self-Org-Ec1}. 
For the first time, FC was used to describe the processes of self-organization, where effects of self-organization with memory is considered, in 2018
\cite{CNSNS2019}, (see also Chapter 19 in \cite[pp. 364--382]{BOOK-DG-2021}).

The emergence and change of memory can leads to the appearance of hierarchy of relaxation times and qualitatively change the behavior of the systems at the remaining other parameters unchanged. This can lead to the self-organization, 
which is generated only by memory. 

We should note that the discrete maps with memory and nonlocality in time, which are derived from fractional differential equations \cite{Tarasov-Zaslavsky,Tarasov-Map1,Tarasov-Map2,Springer2010,Tarasov-Map3,TT-Logistic,Entropy2021,MMAS2021,CSF2021,BOOK-DG-2021}, can demonstrate qualitatively new types of chaotic and regular behavior. 
Similarly, fundamentally new types of self-organization are expected in processes with nonlocality in time and memory.
This direction of future researches is very important for nonlinear fractional dynamics.



\section{Conclusions}
\label{s5}

This paper started by presenting the present day trends of research in FC. Stemming from the state of scientific knowledge, possible directions for further research were pointed. Associated with this analysis some open problems in FC were also identified. The work was organized in three main sub-areas namely, aspects in mathematical analysis, numerical processing, and applications in physics. Given the limitation of space, this work is necessarily parsimonious, but gives a synthetic vision that can somehow guide researchers in the future.




\section*{Declarations}

\paragraph{Funding.}
The authors have not received any external funding in connection with this manuscript.

\paragraph{Conflicts of interest.}
The authors declare that they have no conflict of interest.

\paragraph{Code availability (software application or custom code).}
Not applicable.

\paragraph{Data availability.}
No datasets were generated during the current study.
The data analysed during the study are available from Elsevier's Scopus database but 
restrictions apply to the availability of these data, which were used under license for 
the current study, and so are not publicly available. Data are however available from 
the authors upon reasonable request and with permission of Elsevier B.V.

\paragraph{Authors' contributions.}
All authors contributed to the design and to the writing of the paper.
All authors read and approved the final manuscript.

\end{document}